%% file: ms.tex
\documentclass[final]{l4dc2026}

\usepackage{amsfonts}
\usepackage{dsfont}
\usepackage{enumerate}
\usepackage{soul}
\usepackage{array}
\usepackage{algorithm}
\usepackage{algorithmic}
\usepackage{caption}
\usepackage[table]{xcolor}

\SetCommentSty{mycommfont}

\DeclareOldFontCommand{\rm}{\normalfont\rmfamily}{\mathrm}

\input{header.tex}

\newcommand{\Nx}{n}
\newcommand{\Nu}{m}
\newcommand{\nx}{n_x}
\newcommand{\nuu}{n_u}
\newcommand{\Kor}{K^{o}}
\newcommand{\ustar}{u^{o}}
\newcommand{\opt}{\star}
\newcommand{\Kdrro}{K^\opt}
\newcommand{\phidrro}{\phi^\opt}
\newcommand{\KcalStar}{\Kcal^\opt}
\newcommand{\aff}{v}
\newcommand{\dist}{w}
\newcommand{\xinit}{x_0}
\newcommand{\policy}{\phi}
\newcommand{\Eset}{\mathds{E}}
\newcommand{\cov}{\textrm{cov}}
\newcommand{\muNom}{\widehat{\mu}}
\newcommand{\SigmaNom}{\widehat{\Sigma}}
\newcommand{\DRROC}{distributionally robust regret optimal control }

\newcommand{\Mball}{\Ccal}

\newcommand{\borel}{\Mcal}
\newcommand{\traceB}[2]{\mathop{\rm Tr}#1(#2 #1)}
\newcommand{\Mvar}{\Lambda}
\newcommand{\MvarA}{\Lambda_1}
\newcommand{\MvarB}{\Lambda_2}
\newcommand{\MvarAbar}{\bar{\Lambda}_1}
\newcommand{\MvarBbar}{\bar{\Lambda}_2}
\newcommand{\idx}{i}
\newcommand{\subgrad}{supergradient}
\newcommand{\Subgrad}{Supergradient}
\newcommand{\subgrads}{supergradients}
\newcommand{\subdiff}{superdifferential}

\title[Distributionally Robust Regret Optimal Control]{Distributionally Robust Regret Optimal Control Under Moment-Based Ambiguity Sets}

\usepackage{times}

\coltauthor{\Name{Feras Al~Taha} \Email{foa6@cornell.edu}\\
\Name{Eilyan Bitar} \Email{eyb5@cornell.edu}\\
\addr School of Electrical and Computer Engineering, Cornell University, Ithaca, NY, 14853, USA.}

\begin{document}

\maketitle

\begin{abstract}%
We consider a class of finite-horizon, linear-quadratic stochastic control problems, where the probability distribution governing the noise process is unknown but assumed to belong to an ambiguity set consisting of all distributions whose mean and covariance lie within norm balls centered at given nominal values.
To cope with  this ambiguity, we  design causal affine control policies to minimize the \emph{worst-case expected regret} over all distributions in the  ambiguity set. 
The resulting minimax optimal control problem is shown to admit an equivalent reformulation as a tractable convex program, which can be interpreted as a regularized version of the nominal linear-quadratic stochastic control problem. 
Based on the dual of this convex reformulation, we develop a scalable projected subgradient method for computing optimal controllers to arbitrary accuracy.
Numerical experiments are provided to compare the proposed method with state-of-the-art data-driven control design methods. 
\end{abstract}

\input{Introduction}
\input{Formulation}
\input{Reformulation}
\input{Computation}
\input{Experiments}
\input{Conclusion}

\acks{This work was supported in part by a postgraduate fellowship from the Natural Sciences and Engineering Research Council of Canada, in part by the Cornell Atkinson Center for Sustainability, and in part by the Bezos Earth Fund. The authors thank Dr. Taylan Kargin for a helpful discussion during the initial stages of this research.}

\bibliography{references}

\input{Appendix}

\end{document}

%% file: header.tex
\newtheorem{thm}{Theorem}
\newtheorem{lem}[thm]{Lemma}

\newtheorem{defn}{Definition}

\newtheorem{rem}{Remark}

\newcommand{\argmin}{\mathop{\rm argmin}}

\newcommand{\diag}{\mathop{\mathrm{diag}}}

\newcommand{\trace}[1]{\mathop{\rm Tr}\left(#1\right)}

\newcommand{\Rset}{\mathbb{R}}
\newcommand{\Sset}{\mathbb{S}}

\newcommand{\Acal}{{\cal A}}
\newcommand{\Bcal}{{\cal B}}
\newcommand{\Ccal}{{\mathcal{C}}}

\newcommand{\Kcal}{{\mathcal{K}}}
\newcommand{\Lcal}{{\cal L}}
\newcommand{\Mcal}{{\cal M}}
\newcommand{\Ncal}{{\mathcal{N}}}
\newcommand{\Pcal}{{\cal P}}

\newcommand{\Scal}{{\mathcal{S}}}

\newcounter{l1}
\newcounter{l2}
\newcounter{l3}
\newcommand{\bdotlist}{\begin{list}{$\bullet$}{}}
\newcommand{\bboxlist}{\begin{list}{$\Box$}{}}
\newcommand{\bbboxlist}{\begin{list}{\raisebox{.005in}{{\tiny
$\blacksquare$ \ \ }}}{}}
\newcommand{\bdashlist}{\begin{list}{$-$}{} }
\newcommand{\blist}{\begin{list}{}{} }
\newcommand{\barablist}{\begin{list}{\arabic{l1}}{\usecounter{l1}}}
\newcommand{\balphlist}{\begin{list}{(\alph{l2})}{\usecounter{l2}}}
\newcommand{\bAlphlist}{\begin{list}{\Alph{l2}.}{\usecounter{l2}}}
\newcommand{\bdiamlist}{\begin{list}{$\diamond$}{}}
\newcommand{\bromalist}{\begin{list}{(\roman{l3})}{\usecounter{l3}}}

%% file: Introduction.tex
\section{Introduction} \label{sec:introduction}

A central challenge in the control of dynamic systems subject to uncertain disturbances is the need to make decisions based on incomplete information about the future. In the stochastic control framework, the disturbance is modeled as a stochastic process with a known probability distribution, and a causal control policy is designed to minimize an expected cost criterion under the given distribution. However, stochastic control methods are known to be sensitive to deviations from the assumed distribution, which can significantly degrade performance. In practice, the true distribution is typically unknown due to estimation errors, arising when the distribution is inferred from a finite sample of past disturbances or because of distribution shift  that may occur in nonstationary environments. 

To address these challenges, many distributionally robust approaches to control design have been developed, where the unknown distribution is assumed to belong to a specified set of possible distributions, termed the ambiguity set. Instead of optimizing for a single assumed distribution, these methods seek to enhance the robustness by optimizing against the worst-case distribution within the given ambiguity set.
A variety of ambiguity sets have been studied, including those based on the Wasserstein distance \citep{brouillon2025distributionally,hakobyan2024wasserstein,kim2023distributional,taskesen2024distributionally,yang2020wasserstein}, and those defined through fixed or bounded moments \citep{pmlr-v120-coppens20a, delage2010distributionally, elghaoui2003worst, lee2024performance, li2021distributionally, mark2022recursively, mcallister2024distributionally, van2015distributionally}, to name a few. 
However, the pessimism of distributionally robust control methods can lead to overly cautious control policies and performance loss when the actual disturbance distribution is less severe than anticipated.

To temper the conservatism of (distributionally) robust control design methods with a more optimistic view of the uncertain disturbance process, several recent papers have proposed \emph{worst-case regret minimization} as an alternative to worst-case cost minimization to design control policies for uncertain systems. 
For example, in \citep{sabag2021regret, GoelTAC2023},  the uncertain disturbance process is assumed to lie within a norm-bounded uncertainty set, and causal linear controllers are designed to minimize the worst-case regret across all possible disturbance process realizations within this set.
In this robust regret optimal control framework, \emph{regret} is defined as the excess cost incurred by a causal controller (in response to a particular realization of the disturbance process) relative to the minimum achievable cost by a noncausal controller  (with perfect foresight of the disturbance process at the outset). 
This framework is extended in \citep{liu2024robust} to incorporate uncertainty in the plant dynamics. 

Taking a stochastic view of the disturbance process, \citet{altaha2023distributionally} proposes a distributionally robust approach to regret optimal control, where causal linear controllers are designed to minimize the \textit{worst-case expected regret} over a type-2 Wasserstein ball of disturbance distributions.
Focusing on the finite-horizon setting,  the resulting minimax optimal control problem is shown to admit an equivalent reformulation as a semidefinite program whose size scales polynomially with the control horizon. This framework was further extended to the partially observed setting in \citep{hajar2023wasserstein} and to the infinite-horizon setting in \citep{kargin2024wasserstein, kargininfinite}. 

\textbf{Main Contributions:}
In this paper, we consider the design of distributionally robust regret optimal controllers for the finite-horizon setting  using a novel class of ambiguity sets that consist of distributions with Euclidean norm-bounded means and Schatten norm-bounded covariance matrices.
Focusing on causal affine control policies, we leverage the structure of the proposed ambiguity set to reformulate the minimax optimal control problem as a tractable, finite-dimensional convex program that can be interpreted as a linear-quadratic stochastic control problem with Schatten norm-based regularization of the control policy.
The proposed framework is also  shown to generalize the robust regret optimal control framework studied in \citep{sabag2021regret,GoelTAC2023}.

While the convex reformulation provided in this paper is shown to admit an equivalent representation as a semidefinite program (SDP), solving large-scale SDPs via interior-point methods can be computationally prohibitive in practice.
To address this limitation, we develop a projected subgradient algorithm to efficiently compute controllers to any desired accuracy. 
Finally, we present numerical experiments showing that, when only a finite training sample from the ground-truth distribution is available for control design, the  class of distributionally robust regret optimal controllers proposed in this paper can deliver superior out-of-sample performance compared to other data-driven control design methods in the literature, such as the distributionally robust methods based on Wasserstein ambiguity sets in \citep{altaha2023distributionally}.

\textbf{Paper Organization:} 
The rest of the paper is organized as follows. Sec. \ref{sec:formulation} formulates the distributionally robust regret optimal control problem.
Sec. \ref{sec:reformulation} provides an equivalent reformulation for this problem as a finite-dimensional convex program.
Sec. \ref{sec:opt_methods} presents a  projected subgradient algorithm to efficiently solve this class of problems. 
Sec. \ref{sec:experiments} provides numerical experiments. Sec.~\ref{sec:conclusion} concludes the paper. 
The supplementary material contains omitted proofs (App. \ref{sec:proofs}), SDP reformulation details (App. \ref{sec:sdp}), implementation details of the proposed algorithm (App. \ref{sec:algo_details}), and additional numerical experiments (App.~\ref{sec:additional_exp}). 

%% file: Formulation.tex
\textbf{Notation and Terminology:} 
Let $\Rset$ and $\Rset_+$ denote the set of real numbers and nonnegative real numbers, respectively.  Let $\Sset^n$ denote the set of all symmetric matrices in $\Rset^{n \times n}$. Denote the cone of $n \times n$  real symmetric positive definite (resp. semidefinite) matrices by $\Sset_{++}^{n}$ (resp. $\Sset_+^{n}$). 
Given matrices $A, B \in \Sset^n$, the relation $A \succ B$ (resp. $A\succeq B$) means $A - B \in \Sset^n_{++}$ (resp. $A-B\in\Sset_+^n$). 
Given matrices $A \in \Rset^{m \times p}$ and $B \in \Rset^{n \times p}$,  let $(A, \, B) \in \Rset^{ (m + n) \times p}$ denote the matrix formed by stacking $A$ and $B$ vertically.
Let $\Pi_\Scal(A) := \argmin_{X \in \Scal} \|X-A\| $ denote the orthogonal projection  of the matrix $A\in\Rset^{n \times m}$ onto the set $\Scal\subset \Rset^{n\times m}$ (with respect to the Frobenius inner product).
Let $\borel(\Rset^{n})$ be the set of Borel probability measures on $\Rset^{n}$. 

\section{Problem Formulation} \label{sec:formulation}
Consider a discrete-time, linear time-varying system:
\begin{align} \label{eq:dyn}
    x_{t+1} = A_t x_t + B_t u_t + \dist_t,
\end{align}
which evolves over a finite number of time periods $t=0,\dots,T-1$. Here, $x_t\in\Rset^{\nx}$ is the \textit{system state}, $u_t\in\Rset^{\nuu}$ is the \textit{control input}, and $\dist_t\in \Rset^{\nx}$ is the \textit{disturbance} acting on the system at time $t$. 
We assume that the system matrices $A_t\in\Rset^{\nx\times \nx}$ and $B_t\in\Rset^{\nx\times \nuu}$ are known at the outset and that the system state is perfectly observed.
The initial state $\xinit$ and disturbances $\dist_0,\dots,\dist_{T-1}$ are potentially dependent random variables whose joint distribution $P$ is \emph{unknown}, but assumed to belong to a given compact set of distributions denoted by $\Pcal$ and termed the \textit{ambiguity set}.

We define the associated  \textit{state}, \textit{input}, and \textit{disturbance trajectories} as $x := (\xinit,\,\dots,\,x_T)\in\Rset^{\Nx}$, $u:= (u_0,\,\dots,\,u_{T-1})\in\Rset^{\Nu}$, $\dist := (\xinit,\, \dist_0,\, \dots,\,\dist_{T-1})\in\Rset^{\Nx}$, where $\Nx:=\nx ( T+1 )$ and $\Nu:=\nuu T$.
Note that we have included the initial state $\xinit$ as the first element of the disturbance trajectory $\dist$.
Using these definitions, we can express the dynamics \eqref{eq:dyn} in terms of a causal linear mapping from the input and disturbance trajectories to the state trajectory: 
\begin{align} \label{eq:IO}
    x = Fu + G\dist. 
\end{align}
Here, $F\in\Rset^{\Nx \times \Nu}$ and $G\in\Rset^{\Nx \times \Nx}$ are block lower-triangular (causal) matrices, which are straightforward to construct from the given system matrices $\{A_t, B_t\}_{t=0}^{T-1}$. 

The cost incurred by an input trajectory $u$ and disturbance trajectory $\dist$ is defined as
\begin{align} \label{eq:cost}
    J(u,\,\dist) := x^\top Q x + u^\top R u,
\end{align}
where $Q\in\Sset_+^{\Nx}$ and $R\in\Sset_{++}^{\Nu}$.  We consider  causal affine disturbance feedback controllers
\begin{align} \label{eq:aff} 
    u_t = \aff_t \, + \,  \sum\nolimits_{k=0}^t K_{t,k} \dist_{k-1}, \quad t=0,\dots, T-1,
\end{align}
where $\dist_{-1}:=\xinit$. Here,  $\aff_t \in \Rset^{\nuu}$ denotes the open-loop control input at time $t$, and $K_{t,k} \in \Rset^{\nuu\times\nx}$ denotes the feedback control gain applied to the disturbance $\dist_{k-1}$ at time $t$.
One can also express the control policy \eqref{eq:aff} as an affine mapping $ u = \policy(\dist) :=   K\dist + \aff$ from the disturbance trajectory $w$ to the input trajectory $u$.
Here, $(\aff,\,K) \in \Rset^{\Nu} \times \Lcal$ are the policy parameters, and  $\Lcal\subseteq \Rset^{\Nu\times \Nx}$ represents the space of all block lower triangular matrices that respect the causal disturbance feedback structure specified in \eqref{eq:aff}. Using these definitions, the  set of admissible control policies can be defined as 
\begin{align}
    \Acal := \big\{ \phi(w)  =  Kw + v  \, \big|  \,  (\aff,K) \in \Rset^{\Nu} \times \Lcal  \big\}.
\end{align}

The class of causal affine disturbance feedback control policies is known to be equivalent to the class of causal affine state feedback control policies \citep{goulart2006optimization}.
Specifically, given a disturbance feedback control policy $\phi(w) = Kw + v$, where $\phi \in \Acal$, the change of variables $L:=(I+KG^{-1}F)^{-1}KG^{-1}$ and $c:=(I+KG^{-1}F)^{-1} \aff$ induces a causal state feedback controller $\psi(x) := Lx + c$ that satisfies $\psi \in \Acal$ and $\psi(x) = \phi(w)$ for all  $\dist\in\Rset^{\Nx}$ \citep{lin2019convex, skaf2010design}.
In this paper, we have adopted the disturbance feedback parametrization, because it ensures convexity of the cost function \eqref{eq:cost} with respect to the disturbance feedback controller parameters, whereas the state feedback parameterization results in a nonconvex dependence of the cost function on the state feedback controller parameters. 

\subsection{Distributionally Robust Regret Optimal Control}
This paper is focused  on the design of affine controllers to minimize the worst-case expected regret over all  distributions within a given ambiguity set $\Pcal$.
We define regret as the difference between the cost incurred by a controller and the minimum cost achievable by an optimal noncausal (clairvoyant) controller with perfect knowledge of the disturbance trajectory at the outset.
More specifically,  we define the (ex-post) \textit{regret} incurred by an input trajectory $u$ in response to a disturbance $\dist$ as
\begin{align} \label{eq:reg}
    R(u,\,\dist) := J(u,\,\dist) - J(\ustar(\dist),\,\dist),
\end{align}
where $\ustar:\Rset^{\Nx}\to\Rset^{\Nu}$ denotes the \textit{optimal noncausal controller}, which is defined as $\ustar(\dist) := \argmin_{u \in \Rset^{\Nu}} J(u, \, \dist).$
Building on this definition of regret, we consider the \textit{distributionally robust regret optimal control} problem, as originally proposed in \citep{altaha2023distributionally}: 
\begin{align} \label{eq:mroc}
    \inf_{ \policy \in \Acal} \sup_{P \in  \Pcal}   \Eset_P \left[ R(\policy(\dist),\, \dist) \right].
\end{align}
Compared to worst-case expected cost minimization, the regret-based approach \eqref{eq:mroc} can yield less conservative controllers that deliver more balanced performance across the ambiguity set, although this may come at the expense of a larger worst-case expected cost.

In \citep[Theorem 11.2.1]{hassibi1999indefinite}, the optimal noncausal controller is shown to be both unique and linear in the disturbance trajectory, and is given~by 
\begin{align} \label{eq:Kor}
    \ustar(\dist) = \Kor \dist, 
\end{align}
where $\Kor := -(R + F^\top Q F)^{-1} F^\top Q G.$  Using this explicit characterization of the optimal noncausal controller, the expression for regret given in \eqref{eq:reg} can be simplified to $R(u,\,\dist) = (u - \Kor\dist )^\top D (u- \Kor\dist),$
where $D := R + F^\top QF \succ 0$ \citep{GoelTAC2023, martin2022safe}. It follows that the distributionally robust regret optimal control problem \eqref{eq:mroc} can be reformulated as
\begin{align} \label{eq:mro_reform}
    \inf_{ \policy \in \Acal} \sup_{P \in  \Pcal}   \Eset_P \left[  (\phi(w) - \Kor\dist )^\top D (\phi(w)- \Kor\dist) \right].
\end{align}

\vspace{-.3in}
\subsection{Ambiguity Set} \label{sec:amb_set}

Based on the reformulation provided in \eqref{eq:mro_reform}, the regret incurred by an affine policy $\phi \in \Acal$ is a quadratic function of the disturbance trajectory $w$. 
As a result, the expected regret only depends on the underlying probability distribution governing the disturbance trajectory through its first and second moments. 
With this in mind, we define the ambiguity set as a set of distributions with unknown mean and covariance, which are assumed to lie within a Euclidean norm ball and a Schatten $p$-norm ball, respectively. The Schatten $p$-norm is defined as follows.
\begin{defn}[Schatten $p$-norm] \rm For  $p \in [1, \infty]$, the \textit{Schatten $p$-norm} of  $A \in \Rset^{m \times n}$ is defined as 
\begin{align*}
  \|A\|_{p} : = \begin{cases}    \big( \sum_{i=1}^d \sigma_i^p \big)^{1/p}, &  p \in [1, \infty),\\
    \sigma_1, & p = \infty,
    \end{cases} 
\end{align*}
where $d := \min\{m,n\}$ and $\sigma_1 \geq \sigma_2 \geq \cdots \geq \sigma_d$ denote the singular values of the matrix $A$.
\end{defn}
Note that the Schatten $p$-norm reduces to the \emph{nuclear norm} for $p = 1$, the \emph{Frobenius norm} for $p = 2$, and the \emph{spectral norm} for $p = \infty$.
Additionally, for any $p \in [1, \infty]$, the dual norm associated with the Schatten $p$-norm is the Schatten $q$-norm, where $q \in [1, \infty]$ is such that $\frac{1}{p} + \frac{1}{q} = 1$ \citep{bhatia2013matrix}. 

Given a \emph{nominal mean-covariance pair}  $(\muNom, \SigmaNom) \in \Rset^{\Nx} \times \Sset_+^{\Nx}$,  we define the \emph{ambiguity set} of distributions as \vspace{-0.3cm}
\begin{align} \label{eq:amb_set}
    \Pcal := \big\{ P \in \borel( \Rset^{\Nx}) \ \big| \  \|\mu-\muNom\|_2^2 \leq r_1, \ \|\Sigma-\SigmaNom\|_p  \leq r_2 \big\},
\end{align}
where $\mu := \Eset_P[\dist]$ and $\Sigma := \Eset_P[(\dist-\mu)(\dist-\mu)^\top]$ denote the mean and covariance of the distribution $P$. 
Crucially, this class of ambiguity sets allows for temporally correlated and non-Gaussian disturbance processes, enabling robustness to a broad class of distributions. 

\begin{rem}[Interpretation of Ambiguity Set Parameters] \rm One can interpret the nominal mean $\muNom$ and covariance matrix $\SigmaNom$ as ``point estimates'' of the unknown mean and covariance matrix, and the radii of their respective norm balls, $r_1 \in \Rset_+$ and $r_2 \in \Rset_+$,  as reflecting the degree of ``confidence'' one has in these estimates. 
When $r_1 = r_2 = 0$, problem \eqref{eq:mro_reform} becomes an ``ambiguity-free'' stochastic control problem.
The Schatten-norm order $p$ used to define the covariance ball determines  the shape of the ambiguity set and the form of regularization that arises when solving problem \eqref{eq:mro_reform}. 
When the nominal mean and covariance  $(\muNom, \SigmaNom)$ are learned from data, the parameters $(r_1, \,  r_2, \,  p)$ can be tuned using cross-validation methods to improve controller out-of-sample performance.
\end{rem}

%% file: Reformulation.tex
\section{Convex Reformulation} \label{sec:reformulation}

In this section, we provide an explicit solution to the worst-case expectation problem (the inner maximization) in  \eqref{eq:mro_reform}, which enables the reformulation of problem \eqref{eq:mro_reform} as a tractable convex program.

\subsection{Reformulating the Worst-Case Expected Regret}

The expected regret of  an affine control policy $\phi(w) = Kw + v$ can be expressed as $\Eset_P[R(\policy(\dist),\dist)] = \traceB{\big}{D \, \Eset_P[(\Delta\dist+\aff)(\Delta\dist+\aff)^\top]}$,
where $\Delta := K-\Kor$. Using this reformulation, the worst-case expected regret can be decomposed into two separate maximization problems---one over the disturbance mean $\mu$ and the other over the disturbance covariance $\Sigma$:
\begin{align}
        \label{eq:mroc_separable}  & \max_{\mu\in\Rset^{\Nx} \, : \, \|\mu-\muNom\|_2^2 \le r_1}   (\Delta \mu +\aff)^\top D (\Delta\mu+\aff) \ \ + \max_{ \Sigma\in\Sset_+^{\Nx} \, : \, \|\Sigma-\SigmaNom\|_p \le r_2}  \trace{ \Sigma \Delta^\top D \Delta}.
\end{align}
Note that optimizing over the disturbance mean $\mu$ requires maximizing a convex quadratic function, which results in a nonconvex optimization problem.
While this nonconvexity might lead one to conclude that problem \eqref{eq:mroc} is computationally intractable,  the following lemma  resolves this issue.
\begin{lem} \rm \label{lem:minimax_xy}
    Let $r\ge0$, $y_0\in\Rset^{\Nx}$, $C\in\Rset^{\Nu\times \Nx}$, and $D\in\Sset^{\Nu}_{++}$.
    It holds that 
    \begin{align} \label{eq:minmax_xy}
        \min_{x\in\Rset^{\Nu}} \ \max_{ y \in \Rset^{\Nx} \, : \, \|y-y_0\|_2^2\le r}  (x+Cy)^\top D (x+Cy)  \ =\ r \|C^\top D C\|_\infty,
    \end{align}
    where the optimal value of the minimax problem is attained by the solutions $(x^\opt,\,y^\opt)=(-C y_0,\, y_0 \pm \sqrt{r} \xi)$ and $\xi$ denotes a leading unit eigenvector of the matrix $C^\top D C$. 
\end{lem}
Using  Lemma \ref{lem:minimax_xy} and the decomposition of the worst-case expected regret in \eqref{eq:mroc_separable}, we can now derive an equivalent convex reformulation of the \DRROC problem \eqref{eq:mroc}.

\begin{thm} \label{thm:main_result} \rm 
    Let $r_1,r_2\in\Rset_+$. The following control policy is an optimal solution 
    to problem \eqref{eq:mroc}: 
    \begin{align} 
        \phidrro(\dist) &:= \Kdrro(\dist-\muNom)  \, +  \,  \Kor\muNom,
    \label{eq:phi_star}
    \end{align}
    where $\Kor = -(R + F^\top Q F)^{-1} F^\top Q G$ is the optimal noncausal control policy, and  
    \begin{align} 
    \label{eq:K_star}
        \Kdrro &\!\in \argmin_{K\in\Lcal}  \big\{ \! \traceB{\big}{\SigmaNom C(K)} \, + \, r_1 \|C(K)\|_\infty +  r_2 \|C(K)\|_q \big\},
    \end{align}
where $C(K):=(K-\Kor)^\top D (K-\Kor)$ and $\|\cdot\|_q$ denotes the dual norm of the Schatten $p$-norm.
\end{thm}

The \DRROC policy $\phidrro(\dist) = \Kdrro(\dist-\muNom) + \Kor\muNom$  may be interpreted as follows. The \emph{open-loop term} $\Kor \muNom$, which can be precomputed offline, equals the expected action of the optimal noncausal control policy \eqref{eq:Kor} under the nominal disturbance mean $\muNom$. 
The causal \emph{feedback term} $\Kdrro(\dist-\muNom)$  compensates online for any deviations between the realized disturbance and the nominal disturbance mean. 

\subsection{Controller Regularization}  \label{sec:regularization} 

The reformulation  in \eqref{eq:K_star} can be interpreted as a \emph{regularized} linear-quadratic stochastic control problem. 
The first term in the  objective is the expected regret of a linear controller $K$ in response to a zero-mean disturbance $w$ with covariance $\SigmaNom$. 
By adding the expected cost of the optimal \emph{noncausal} controller under the same disturbance distribution,  the optimization in \eqref{eq:K_star}  can be rewritten as $    \Kdrro \! \in \! \argmin_{K\in\Lcal}  \Eset_{\bar{P}}[J(Kw,w)]  +  r_1 \|C(K)\|_\infty \!+  r_2 \|C(K)\|_q $, where $\bar{P} \in \borel(\Rset^{n})$ is any distribution satisfying $\Eset_{\bar{P}}[w] \!=\! 0$  and  $\Eset_{\bar{P}}[ww^\top] \!= \!\SigmaNom$. 
The leading term represents the nominal expected cost of the controller $K$ under the distribution $\bar{P}$. 
The second and third terms penalize the controller's sensitivity to distributional uncertainty in the disturbance mean and covariance, respectively. 

The spectral norm penalty $r_1 \|C(K)\|_\infty$ arises from  mean uncertainty and coincides with  the worst-case regret incurred by the controller $K$ over all  norm-bounded disturbances. Specifically, 
\begin{align*}
     \max_{\|\bar{w}\|_2\le 1} R(K \bar{w}, \bar{w}) = \max_{\|\bar{w}\|_2\le 1} \bar{w}^\top C(K) \bar{w} = \|C(K)\|_\infty,
\end{align*}
so minimizing this term with respect to $K$ corresponds  to solving the robust regret optimal control problem $ \min_{K \in \Lcal} \max_{\|\bar{w}\|_2\le 1} R(K \bar{w}, \bar{w})$, as studied in \citep{GoelTAC2023,sabag2021regret}.  

The Schatten $q$-norm penalty $r_2 \|C(K)\|_q $ stems from uncertainty in the  covariance matrix, which is assumed to lie within a Schatten $p$-norm ball, with  $1/p + 1/q = 1$.
Thus, when the covariance uncertainty set is defined via  the nuclear norm ($p=1$),  the dual Schatten $q$-norm reduces to the spectral norm, and the two regularization terms in \eqref{eq:K_star} combine to form a single penalty of the form $(r_1 + r_2) \|C(K)\|_\infty$.
Conversely, when the covariance matrix is bounded in spectral norm ($p = \infty$), the dual Schatten $q$-norm  reduces to the nuclear norm, yielding a regularization term of the form $ r_2 \|C(K)\|_1 = r_2 \trace{C(K)}$,  which coincides with the expected regret of the linear controller $K$ under a zero-mean disturbance process with isotropic covariance $r_2 I$. 
When the cost function \eqref{eq:cost} is separable across time periods (i.e., $Q$ and $R$ are block-diagonal),  the quantity $ \trace{C(K)}$ is minimized by the classical linear-quadratic regulator (see App. \ref{sec:lqr} of the supplementary material). 

With these insights, the  \DRROC problem can be understood as \emph{interpolating} between three canonical control design objectives: (i) the nominal linear-quadratic stochastic control problem under the  distribution $\bar{P}$ (as $r_1, r_2 \rightarrow 0$);  (ii) the robust regret optimal control problem (as $r_1 \rightarrow \infty)$; and (iii)  the classical linear-quadratic regulator (as $r_2 \rightarrow \infty$ and $p = \infty$). Thus, by tuning the ambiguity set parameters $(r_1,r_2,p)$, one can design controllers that trade off nominal performance and robustness to different forms of distributional shift.

%% file: Computation.tex
\section{Dual Projected \Subgrad{} Method} \label{sec:opt_methods}

The convex program provided in Theorem \ref{thm:main_result} can be recast as a semidefinite program (SDP) whose size grows polynomially with $\nx$, $\nuu$, and $T$ (see App. \ref{sec:sdp}).
While  SDPs can be solved in polynomial time using interior point methods, their memory and computational requirements  grow rapidly with the problem size, limiting their scalability to large~systems.
To address this, we develop a projected \subgrad{} method to solve the dual of problem \eqref{eq:K_star}. The algorithm reduces to projected gradient ascent  when the nominal covariance  $\SigmaNom$ is positive definite. Numerical results in App. \ref{sec:additional_exp} show that the proposed method delivers substantial speedups  over a state-of-the-art interior point solver.

\subsection{Dual Reformulation} 
To derive the dual of problem \eqref{eq:K_star}, it is helpful to define the following  Cartesian product of covariance norm balls:
$\Ccal \!\! := \!\! \big\{ \Sigma\in\Sset_+^{\Nx} \,\big\vert\, \|\Sigma\|_1  \!\le \! r_1 \big\} \!\! \times \!\! \big\{ \Sigma\in\Sset_+^{\Nx} \,\big\vert\, \|\Sigma-\SigmaNom\|_p  \!\le \! r_2 \big\}.$
With this definition, problem \eqref{eq:K_star} can be expressed as the minimax problem $\min_{K \in \Lcal}  \, f(K),$ where $f(K) :=  \max_{\Mvar\in\Mball} \traceB{\big}{ G(K)^\top \Mvar }$ and $G(K)\!:=\!(C(K),C(K))$.  
Using Sion's minimax theorem \citep{sion1958general}, the order of minimization and maximization can be reversed, leading to the dual reformulation:  
\begin{align} \label{eq:max_sigma}
    \min_{K \in \Lcal}  \, f(K) = \max_{\Mvar\in\Mball} \, g(\Mvar), 
\end{align}  
where $g(\Mvar) :=  \min_{K \in \Lcal} \traceB{\big}{ G(K)^\top \Mvar}$.  
In the following theorem, we establish several properties of the dual function $g$ that  will facilitate the design of a dual projected \subgrad{} method. 
\begin{thm} \label{cor:max_sigma} \rm
The  dual function $g$ has the following properties:
\begin{enumerate}[(i)] 
\setlength{\itemsep}{5pt} 
\setlength{\parskip}{0pt}  
\setlength{\parsep}{0pt}   
    \item \emph{Concavity:} The function  $g$ is concave over $\Sset_+^{\Nx}\times\Sset_+^{\Nx}$.
    
    \item \emph{Superdifferential:}  The \subdiff{} of $g$ at any point $\Mvar\in\Sset_+^\Nx\times\Sset_+^\Nx$ is nonempty and contains the following set of \subgrads: 
    \begin{align} \label{eq:subdiff} 
        \hspace{-0.5cm}\partial g(\Mvar) \supseteq {\rm conv}\big(\big\{ G(K) \mid K \in \KcalStar(\Mvar)  \big\}\big),
    \end{align}
    where $\KcalStar(\Mvar) := \argmin_{K \in \Lcal} \traceB{\big}{G(K)^\top \Mvar}$ is the nonempty set of minimizers of $ \traceB{\big}{G(\cdot)^\top \Mvar} $.

    \item \emph{Differentiability:} The function  $g$ is differentiable at any  $\Mvar:=(\MvarA, \MvarB) \in \Sset_+^{\Nx}\times\Sset_+^{\Nx}$ satisfying $\MvarA+\MvarB \succ 0$.
\end{enumerate}
\end{thm}

\subsection{Dual Projected \Subgrad{} Method}
 Using Theorem \ref{cor:max_sigma}, we  provide a projected \subgrad{} method to solve the dual  problem in \eqref{eq:max_sigma}. Given the initial condition $ \Mvar^0 := (0,\,\SigmaNom)$ and a sequence of positive step sizes $\{\eta^\idx\}_{\idx=0}^\infty$, compute
\begingroup
\allowdisplaybreaks
\begin{align} 
    \label{eq:subgrad_update_1}  K^\idx \, &\in \argmin_{K\in\Lcal} \traceB{\big}{G(K)^\top \Mvar^\idx},\\
    \label{eq:subgrad_update_2} \Mvar^{\idx+1} &= \Pi_{\Mball} \big(\Mvar^\idx + \eta^\idx G(K^\idx) \big),
\end{align}
\endgroup
until a stopping criterion is satisfied.  
Given a nonsummable, diminishing sequence of step sizes and bounded \subgrads, the sequence of iterates $\{\Mvar^\idx\}_{\idx=0}^\infty$
generated by the dual projected \subgrad{} method \eqref{eq:subgrad_update_1}-\eqref{eq:subgrad_update_2} is guaranteed to  converge asymptotically to an optimal solution of the dual problem \eqref{eq:max_sigma} \citep[Proposition 3.2.6]{bertsekas2015convex}.
Similar convergence guarantees hold for certain classes of adaptive step size rules \citep{bertsekas2015convex}.
Implementation details for this method, including supergradient computation, projection,  and stopping criterion, are provided in App. \ref{sec:algo_details}.

\subsection{Reduction to Dual Projected Gradient Ascent} \label{sec:pga}
The following theorem shows that the proposed projected \subgrad{} method \eqref{eq:subgrad_update_1}-\eqref{eq:subgrad_update_2} reduces to a projected gradient ascent method when the nominal covariance matrix $\SigmaNom$ is positive definite. 

\begin{thm} \rm   \label{thm:pd_iterates}
   The sequence of iterates $\{\Mvar^\idx = (\MvarA^\idx,\,\MvarB^\idx)\}_{i=0}^\infty$ generated by the projected \subgrad{} method \eqref{eq:subgrad_update_1}-\eqref{eq:subgrad_update_2} satisfies $\MvarB^\idx \succeq \SigmaNom$ for all $i\ge 0$.
   Moreover, if $\SigmaNom\succ 0$, then the dual function $g$ is differentiable at every iterate in this sequence.
\end{thm}

By Theorem \ref{thm:pd_iterates}, if the dual projected \subgrad{} method is initialized at  $\Mvar^0=(0,\SigmaNom)$ with a positive definite nominal covariance matrix $\SigmaNom$, then the dual function $g$ is differentiable at all subsequent iterates $\{\Mvar^i\}_{i=0}^\infty$.
Under this condition, the projected \subgrad{} method reduces to a projected \emph{gradient ascent} method. This reduction raises the possibility of using iterative optimization methods that leverage smoothness to accelerate convergence, including momentum-based updates, line search, and quasi-Newton methods \citep{Bertsekas99}.
\vspace{-0.05cm}

%% file: Experiments.tex
\section{Numerical Experiments} \label{sec:experiments}

We study a damped double integrator system affected by stochastic disturbances over a time horizon of $T=10$, with  cost matrices $Q=I$ and $R=10 I$, and linear time-invariant dynamics $A_t={\small \begin{bmatrix} 1 & 1\\ 0 & 0.05\end{bmatrix}}$ and $B_t={\small \begin{bmatrix} 0\\ 1\end{bmatrix}}$.
The disturbance follow an AR(1) process  $w_t = \rho w_{t-1} + \epsilon_{t}$, where  $w_{-1} := x_0 \sim \Ncal(0,\,I)$ and  $\epsilon_{t}\sim\Ncal(0,\,(1-\rho^2)I)$ for $t=0,\dots,T-1$ are mutually independent random variables. This construction induces an autocovariance $\cov(w_{s-1},\,w_{t-1}) = \rho^{|s-t|} I$, allowing us to control temporal dependence through the parameter $\rho\in[-1,\,1]$.
The true disturbance distribution is assumed to be unknown. Instead, we are given a finite training sample of $N=\Nx+1$\, i.i.d. disturbance trajectories $w^{(1)},\dots,w^{(N)} \in\Rset^{\Nx}$, from which we construct the nominal covariance as $\SigmaNom = \frac{1}{N}  \sum_{i=1}^N w^{(i)} {w^{(i)}}^\top$. We assume no uncertainty in the disturbance mean (i.e., $\muNom = 0$ and $r_1=0$) and let $r$ denote the covariance uncertainty radius. Using these data, we synthesize distributionally robust regret optimal controllers and compare their out-of-sample performance against alternative data-driven baselines. 
Additional experiments are  provided in App. \ref{sec:additional_exp} to examine controller performance in the presence of uncertainty in the disturbance mean and to demonstrate the convergence and scalability of the dual projected \subgrad{} method \eqref{eq:subgrad_update_1}-\eqref{eq:subgrad_update_2}.

\subsection{Proposed and Baseline Control Methods Under Comparison} \label{sec:controllers}
We compare the proposed framework against several representative baselines spanning different ambiguity models and optimization criteria:
\begin{enumerate}
\setlength{\itemsep}{8pt}  
\setlength{\parskip}{0pt}  
\setlength{\parsep}{0pt}   

    \item \textbf{Schatten-norm covariance uncertainty sets (this paper)}:  Controllers that minimize worst-case expected regret over covariance uncertainty sets defined by Schatten $p$-norm bounds. We evaluate three variants corresponding to $p \in \{1, 2, \infty\}$, referred to as \textsc{Nuc-Regret}, \textsc{Frob-Regret}, and \textsc{Spec-Regret}, respectively. 
    
    \item \textbf{Wasserstein ambiguity sets}:  Controllers optimized over type-2 Wasserstein balls centered at a nominal distribution with mean-covariance $(\muNom =0, \SigmaNom)$, based on the methods in  \citep{altaha2023distributionally}. We evaluate both worst-case regret-minimizing and worst-case cost-minimizing controllers, denoted \textsc{Wass-Regret} and \textsc{Wass-Cost}, respectively. The Wasserstein radius is set  to $\sqrt{r}$ so that its units match those of the Schatten-norm covariance radius $r$.
    
    \item \textbf{Sample Average Approximation (SAA):} The controller that minimizes the empirical expected cost under the nominal mean-covariance pair  $(\muNom =0, \SigmaNom)$,  referred to as SAA.
    
    \item \textbf{Optimal Causal Controller:} The optimal causal affine controller that minimizes the expected cost under the true mean-covariance pair \((\mu=0,\Sigma)\), referred to as \textsc{Opt-Causal}.
\end{enumerate}

\begin{figure}[h]
    \centering
    \includegraphics[width=\linewidth,trim=0cm 0.835cm 1.0cm 0.25cm,clip]{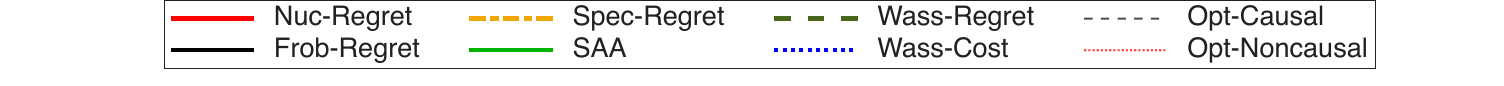}
    \subfigure[Cost versus radius $r$ \label{fig:Ecost_vs_radius}]{
     \includegraphics[width=0.45\linewidth,trim=0.585cm 0.54cm 0.5cm 0.3cm,clip]{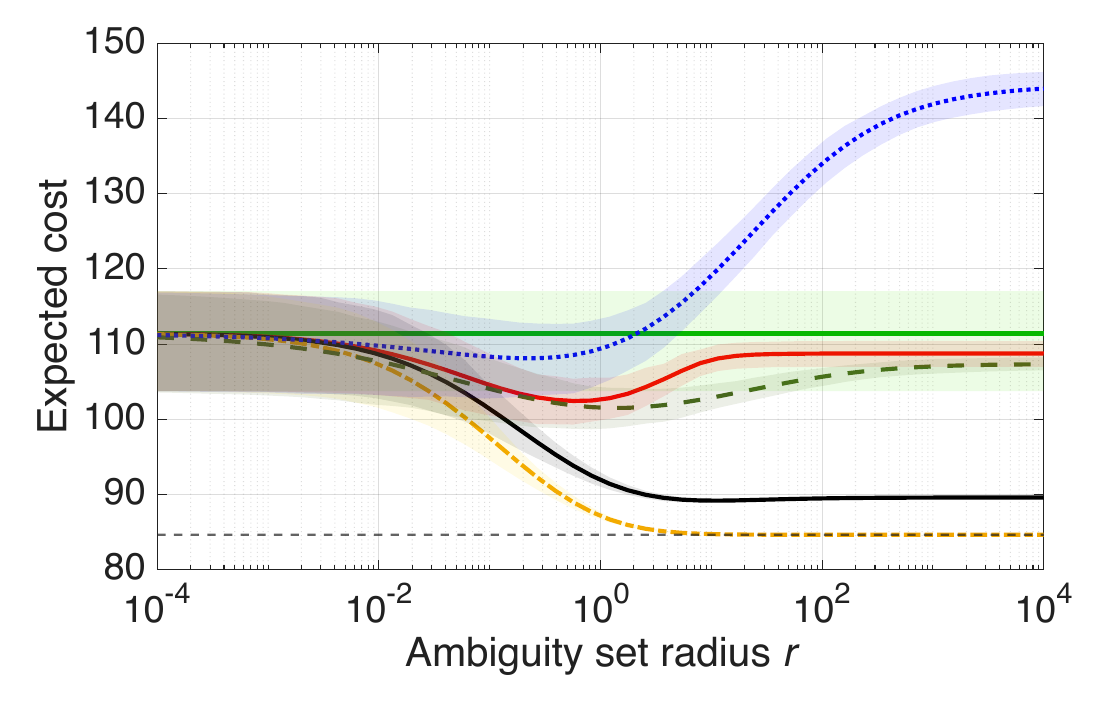}}\hspace{1.0cm}
    \subfigure[Cost versus $\rho$ \label{fig:Ecost_vs_rho}]{
    \includegraphics[width=0.45\linewidth,trim=0.585cm 0.54cm 0.5cm 0.3cm,clip]{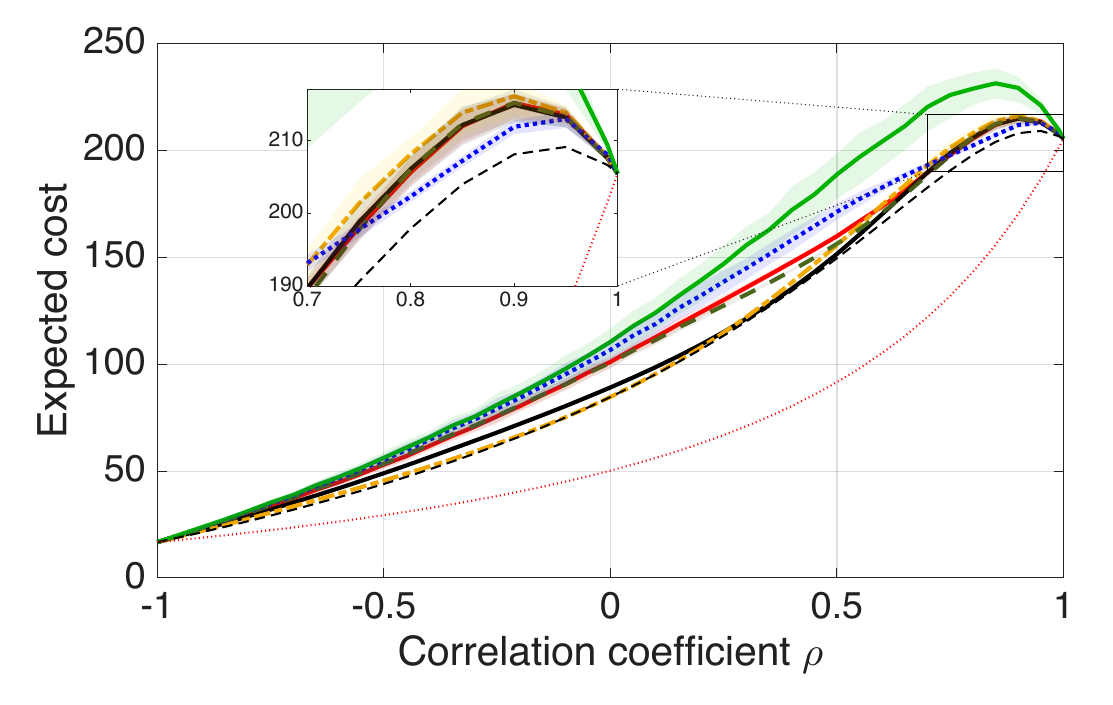}}
    \subfigure[Ex-ante regret versus $\rho$ \label{fig:Eregret_exante_vs_rho}]{
    \includegraphics[width=0.45\linewidth,trim=0.585cm 0.54cm 0.5cm 0.3cm,clip]{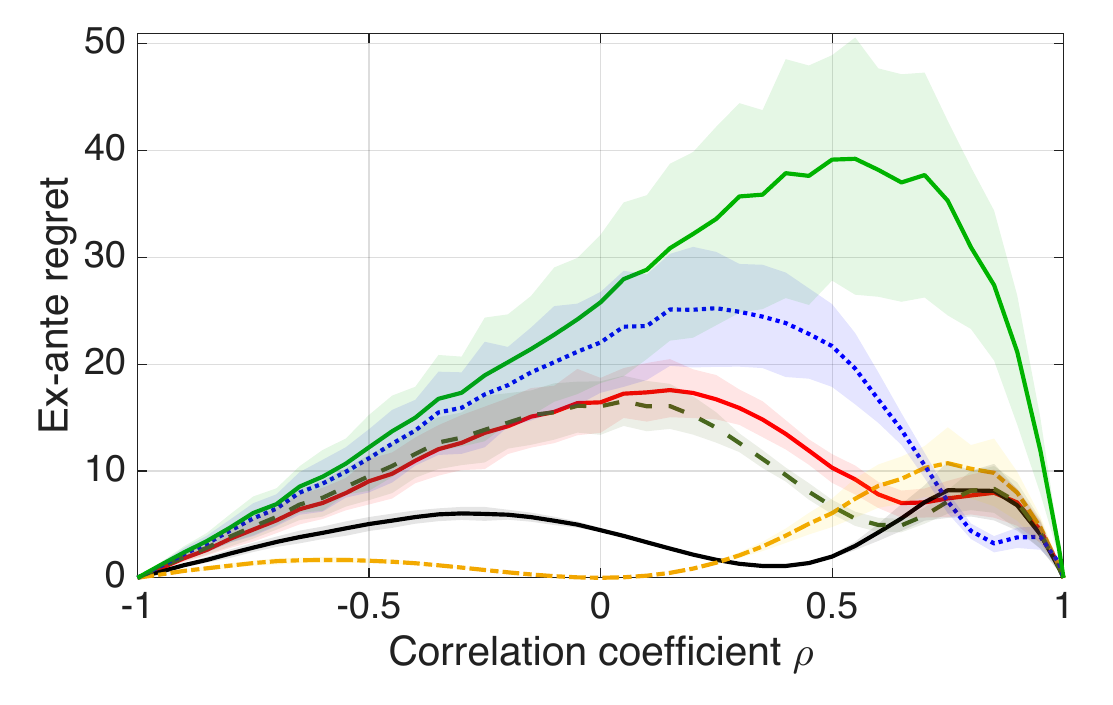} }\hspace{1.0cm}
    \subfigure[Ex-post regret versus $\rho$ \label{fig:Eregret_expost_vs_rho}]{
    \includegraphics[width=0.45\linewidth,trim=0.585cm 0.54cm 0.5cm 0.3cm,clip]{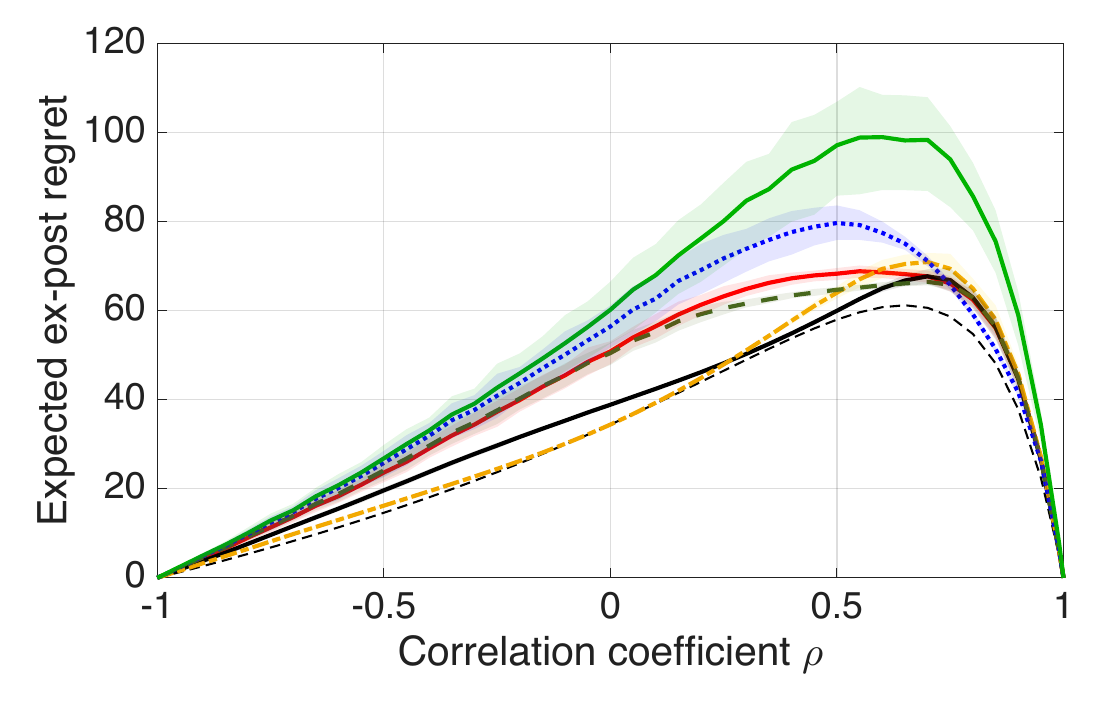} }
    \caption{(a) Expected cost versus the ambiguity radius $r$ for $\rho=0$. (b)-(d) Expected cost, ex-ante regret, and expected ex-post regret versus the correlation parameter $\rho$, where each controller is evaluated at its best performing radius $r$. In the legend, \textsc{Opt-Noncausal} refers to the optimal noncausal controller \eqref{eq:Kor}. All expectations are taken under the ground-truth distribution.  Reported values are  averaged over 200  independent trials, where in each trial the nominal covariance is estimated from a fresh training sample. 
    Shaded regions show the empirical 20\textsuperscript{th}--80\textsuperscript{th} percentile range. 
    }
\end{figure}

\subsection{Out-of-Sample Performance}
We evaluate the out-of-sample performance of all controllers across a family of ground-truth distributions, parameterized by $\rho \in [-1, \, 1]$.
For each value of $\rho$, we conduct 200 independent trials. In every trial, we recompute the controllers using a fresh training sample drawn from the ground-truth distribution.  For ambiguity-set-based methods, we sweep over $r \in [10^{-4}, \, 10^{4}]$, and recompute the controllers for each value of $r$ in this range.

We begin by examining how the ambiguity set radius affects the out-of-sample performance of the different controllers when the disturbance process is  temporally uncorrelated ($\rho=0$). Figure \ref{fig:Ecost_vs_radius} depicts the expected cost of each controller under the ground-truth distribution (averaged over 200 trials) as a function of the ambiguity set radius $r$. When $r=0$, all of the distributionally robust controllers coincide with the SAA controller. As the ambiguity set radius increases, the out-of-sample performance of each method improves relative to the SAA controller, benefiting from distributional robustness. However, as the radius grows large, several controllers become overly conservative and their out-of-sample performance begins to degrade. In particular, for sufficiently large radii, the \textsc{Wass-Cost} controller performs worse than the SAA controller. When comparing the different controllers at their corresponding \emph{optimal radii} (the value of $r$ that yields the lowest out-of-sample expected cost), the \textsc{Spec-Regret} controller outperforms all other methods, achieving an expected cost that matches that of the \textsc{Opt-Causal} controller. This aligns with the claim in Sec. \ref{sec:regularization}  that, as $r\to\infty$, the \textsc{Spec-Regret} controller converges to the classical linear-quadratic regulator, which coincides with the \textsc{Opt-Causal} controller  in this case, because the disturbance process is uncorrelated across time $(\rho = 0)$ and the cost function is separable across time periods.
While these observations highlight the importance of radius selection in shaping out-of-sample performance, the optimal radius cannot be computed in practice, since the ground-truth distribution is unknown. The optimal radius can, however, be approximated using data-driven methods that estimate the out-of-sample expected cost via holdout or cross-validation techniques, as described in \citep[Sec. 7.2.2]{mohajerin2018data}.

Next, we compare the out-of-sample performance of the different controllers across the full range of correlation parameters $\rho \in [-1, \,1]$. For each value of $\rho$, we calculate the optimal radius for each distributionally robust controller over the range $r \in [10^{-4}, \, 10^{4}]$. Figures \ref{fig:Ecost_vs_rho}-\ref{fig:Eregret_expost_vs_rho} report the expected cost, ex-ante regret,\footnote{The \emph{ex-ante regret} of a controller is defined as the difference between its expected cost under the ground-truth distribution and the expected cost achieved by the \textsc{Opt-Causal} controller, which knows the ground-truth distribution.} and expected ex-post regret of the resulting controllers as a function of $\rho$ (averaged over 200 trials). 
As seen in Figure \ref{fig:Ecost_vs_rho}, the expected cost of all controllers increases with the correlation parameter $\rho$ before tapering off slightly at very high correlation levels. As $\rho$ increases, disturbances are more likely to accumulate over time, increasing the control effort required for disturbance rejection. However, as $\rho \rightarrow 1$, the disturbance process becomes more predictable, yielding a small reduction in expected cost for $\rho$ close to one. 
When $\rho=\pm1$, the disturbance sequence $\{w_t\}_{t=0}^{T-1}$ is completely determined by the initial state $w_{-1}=x_0$. As a result, the full disturbance trajectory is known from the outset, and every  controller achieves zero ex-post regret.
For the ``least-favorable'' values of $\rho$, the \textsc{Wass-Cost} controller performs best, achieving the lowest out-of-sample expected cost over a narrow band of high-correlation parameters. Outside that region, however, the regret-minimizing controllers perform substantially better. Among these, the \textsc{Spec-Regret} and \textsc{Frob-Regret} controllers achieve the best robustness-performance tradeoff across the family of distributions considered, delivering consistently low regret while avoiding the  over-conservatism exhibited by \textsc{Wass-Cost}, as shown in Figures~\ref{fig:Eregret_exante_vs_rho}-\ref{fig:Eregret_expost_vs_rho}.

%% file: Conclusion.tex
\section{Conclusion} \label{sec:conclusion}
In this paper, we provide a framework for synthesizing  distributionally robust regret optimal controllers for linear systems affected by additive stochastic disturbances with unknown but bounded first and second moments. A tractable convex reformulation of the \DRROC problem is provided, along with a projected subgradient algorithm to enable the efficient computation of control policies.

%% file: Appendix.tex
\newpage
\vspace*{.00in}
\begin{center}
    {\LARGE \textbf{Supplementary Material}}
\end{center}

The supplementary material of this paper is organized as follows.
App. \ref{sec:proofs} provides all omitted proofs for the theoretical results in this paper.
App. \ref{sec:sdp} shows how to reformulate the distributionally robust regret optimal control problem \eqref{eq:mroc} as an equivalent semidefinite program.
App. \ref{sec:algo_details} describes implementation details of the proposed dual projected supergradient algorithm.
App. \ref{sec:additional_exp} provides additional numerical experiments. 

\appendix

\section{Deferred Proofs} 
\label{sec:proofs}

\subsection{Proof of Lemma \ref{lem:minimax_xy}}

It is straightforward to verify that the solutions $(x^\opt,\,y^\opt)=(-C y_0,\, y_0 \pm \sqrt{r} \xi)$ are feasible for the minimax problem and attain the value $r \|C^\top D C\|_\infty$.
To complete the proof, it suffices to show that the optimal value of the minimax problem is  upper and lower bounded by $r  \|C^\top D C\|_\infty$. 
Taking $x=-Cy_0$, we can upper bound the optimal value of the outer minimization as 
\begin{align*}
    \min_{x\in\Rset^{\Nu}} \max_{\substack{y\in\Rset^{\Nx}\\ \|y-y_0\|_2^2\le r}} (x+Cy)^\top D (x+Cy)
    &\le \max_{\substack{y\in\Rset^{\Nx}\\ \|y-y_0\|_2^2 \le r}} (y-y_0)^\top C^\top D C (y-y_0)\\
    & \le \max_{\substack{y\in\Rset^{\Nx}\\ \|y-y_0\|_2^2 \le r}} \| C^\top D C\|_\infty \|y-y_0\|_2^2 \\
    &=r \|C^\top D C\|_\infty,
\end{align*}
where the second inequality follows from the Cauchy-Schwarz inequality.

To show that $r \|C^\top D C\|_\infty$ is also a lower bound on the value of the minimax problem, take $y=y_0+\alpha(x) \xi$, where $\alpha(x) := \sqrt{r} \, \textrm{sgn}( (x+C y_0)^\top D C \xi )$ and $\textrm{sgn}(\cdot)$ denotes the sign function with the convention $\textrm{sgn}(0)=1$.
It follows that $y=y_0+\alpha(x) \xi$ is feasible since $\|y-y_0\|_2^2=\|\alpha(x) \xi\|_2^2=\alpha(x)^2 = r$ for all $x\in\Rset^\Nu$.
Hence, the value of the minimax problem can be lower bounded as
\begin{align}
    \nonumber \min_{x\in\Rset^{\Nu}} &\max_{\substack{y\in\Rset^{\Nx}\\ \|y-y_0\|_2^2\le r}} (x+Cy)^\top D (x+Cy)\\
    \nonumber &\ge \min_{x\in\Rset^{\Nu}} \big(x+Cy_0 + \alpha(x) C \xi\big)^\top D \big(x+Cy_0 + \alpha(x) C \xi\big)\\
    \label{eq:lower_min} &= \min_{x\in\Rset^{\Nu}} (x+Cy_0)^\top D (x+Cy_0) +2 \alpha(x) (x+Cy_0)^\top D C \xi + r \|C^\top D C\|_\infty.
\end{align}
Note that the first two terms of the objective of the minimization problem in \eqref{eq:lower_min} are nonnegative both when $(x+C y_0)^\top D C \xi \ge 0$ (which implies $\alpha(x)=\sqrt{r}$), or when $(x+C y_0)^\top D C \xi < 0$ (which implies $\alpha(x)=-\sqrt{r}$).
Since these two terms equal zero when evaluated at the feasible solution $x^\opt=-Cy_0$, it follows that $x^\opt$ is an optimal solution of the minimization problem in \eqref{eq:lower_min}, yielding an optimal value of $r \|C^\top D C\|_\infty$, which completes the proof.

\subsection{Proof of Theorem \ref{thm:main_result}} \label{app:proof_main}

\begingroup
\allowdisplaybreaks
    To prove the desired result, we first use the fact that the objective function in \eqref{eq:mro_reform} is separable in $\Sigma$ and $(\aff,\mu)$ to derive a closed-form solution for the optimal open-loop control input $\aff$ in terms of the feedback control gain matrix $K$. We then solve the inner maximization with respect to $\Sigma$.

    Using the reformulation of the worst-case expected regret provided in \eqref{eq:mroc_separable}, we can rewrite the \DRROC problem \eqref{eq:mroc} as 
    \begin{align} \label{eq:mroc_separable_proof}
        &\min_{K\in\Lcal} \Bigg \{  \min_{\aff\in\Rset^{\Nu}} \, \max_{\substack{\mu\in\Rset^{\Nx} \, : \, \|\mu-\muNom\|_2^2\le r_1}}  (\Delta \mu +\aff)^\top D (\Delta\mu+\aff) \ + \max_{\substack{\Sigma\in\Sset_+^{\Nx} \, : \, \|\Sigma-\SigmaNom\|_p\le r_2}}  \trace{ \Sigma C(K)}  \Bigg \}. 
    \end{align}
    By Lemma \ref{lem:minimax_xy}, the inner minimax problem in $(\aff,\mu)$ has an optimal value of $r_1 \|C(K)\|_\infty$, which is attained by the solution $\aff^\opt=-\Delta \muNom$. 
    To solve the remaining maximization in $\Sigma$, take any fixed $K\in\Lcal$.
    Then, we can upper bound the optimal value of the maximization problem in $\Sigma$ by dropping the positive semidefiniteness condition on $\Sigma$ as follows:
    \begin{align}
       \nonumber \max_{\substack{\Sigma\in\Sset_+^{\Nx} \, : \, \|\Sigma-\SigmaNom\|_p\le r_2}}  \trace{ \Sigma C(K)}
        &\le  \max_{\substack{\Sigma\in\Sset^{\Nx} \, : \, \|\Sigma-\SigmaNom\|_p\le r_2}}  \traceB{\big}{\SigmaNom C(K)} +\traceB{\big}{ (\Sigma-\SigmaNom) C(K) }\\
        &=  \traceB{\big}{\SigmaNom C(K)} + r_2 \|C(K)\|_q, \label{eq:upper1}
    \end{align}
    where the above equality follows from the definition of the dual norm.
    Now, letting $\xi\in\Rset^{\Nx}$ denote a leading unit eigenvector of $C(K)$, define
        \begin{align*}  
        \Sigma^\opt := \begin{cases}
             \displaystyle \SigmaNom + r_2 \|C(K)\|_q^{-(q-1)} C(K)^{q-1}, &  p \in (1, \infty],\\
             \SigmaNom + r_2 \xi\xi^\top, & p=1,
        \end{cases}
    \end{align*}
    where we use the convention $C(K)^{q-1} = I$ when $q=1$.
    It is straightforward to verify that, for each $p\in[1,\infty]$, the covariance matrix $\Sigma^\opt$ attains the upper bound in \eqref{eq:upper1} and is feasible for the original maximization problem in $\Sigma$ since it can be shown to satisfy $\|\Sigma^\opt-\SigmaNom\|_p = r_2$ and $ \Sigma^\opt \in\Sset^{\Nx}_+$. 
    Hence, the upper bound in \eqref{eq:upper1} equals the optimal value of the maximization problem in $\Sigma$, and \eqref{eq:mroc_separable_proof} reduces to the minimization problem in \eqref{eq:K_star}.
    Finally, \eqref{eq:phi_star} follows from substituting in $\aff^\opt=-(\Kdrro-\Kor) \muNom$ and $\Kdrro$ into $\phi(\dist)=K\dist+\aff$, concluding the proof.
\endgroup

\subsection{Reduction to Classical Linear-Quadratic Regulator} \label{sec:lqr}
When the cost function \eqref{eq:cost} is additively separable across time periods, the function $\trace{C(K)}$ is minimized by the classical linear-quadratic regulator (LQR).
To see why this is true,
suppose that $Q$ and $R$ are block-diagonal  with  blocks $Q_t\in\Sset_+^{\nx}$ ($t=0,\dots,T$) and $ R_t\in\Sset_{++}^{\nuu}$ ($t=0,\dots,T-1$), respectively.  Additionally,  let $P\in\borel(\Rset^\Nx)$ be any distribution with  mean $\Eset_P[w]=0$ and   covariance $\Eset_P[ww^\top]=I$. Then, for any $K\in\Lcal$, 
\begin{align}
    \nonumber \trace{C(K)}  = \traceB{\big}{C(K) \Eset_P[ww^\top]} =  \Eset_{P}[w^\top (K-\Kor)^\top D (K-\Kor) w] = \Eset_{P}[R(Kw,\,w) ],
\end{align}
where the final equality follows from the equivalence between \eqref{eq:mroc} and \eqref{eq:mro_reform}. Using  the time-separable cost structure and writing $u = Kw$, this becomes
\begin{align}
    \label{eq:lqr_cost} \trace{C(K)}  =  \Eset_P\Big[x_T^\top Q_T x_T + \sum_{t=0}^{T-1} x_t^\top Q_t x_t + u_t^\top R_t u_t  \Big]   - \Eset_{P}[ J(\ustar(w),\,w)].
\end{align}
 Now, recall that every (strictly) causal linear disturbance feedback policy $u = Kw$ admits an equivalent causal linear  state feedback representation $u = Lx$ for some $L \in \Lcal$ (see Sec. \ref{sec:formulation}).  Since the second term in \eqref{eq:lqr_cost} is independent of $K$, minimizing $\trace{C(K)}$ over $K \in \Lcal$ is equivalent to minimizing $$\Eset_P\Big[x_T^\top Q_T x_T + \sum_{t=0}^{T-1} x_t^\top Q_t x_t + u_t^\top R_t u_t  \Big] $$ over all causal linear state feedback policies $u = Lx$ with $L \in \Lcal$. This is precisely the finite-horizon, perfectly observed linear-quadratic control problem, whose solution is the classical LQR policy \citep{bertsekas2012dynamic}. 

\subsection{Proof of Theorem \ref{cor:max_sigma}}

    \noindent \emph{Proof of (i):} Concavity of the function $g$ on $\Sset_+^{\Nx}\times\Sset_+^{\Nx}$ follows from the fact that $g$ is the pointwise minimum of linear functions.\\

    \noindent \emph{Proof of (ii):}
    Let $\bar{\Mvar} := (\MvarAbar,\,\MvarBbar) \in \Sset_+^{\Nx}\times\Sset_+^{\Nx}$.
    The function $\traceB{\big}{G(K)^\top \bar{\Mvar}} = \traceB{\big}{ (\MvarAbar+\MvarBbar)C(K)}$ is a convex quadratic function in $K$ that is bounded from below by zero.
    This guarantees that the set of minimizers $\KcalStar(\bar{\Mvar}) := \argmin_{K \in \Lcal} \traceB{\big}{G(K)^\top \bar{\Mvar}}$ is nonempty.
    Hence, the set defined in the right-hand side of the inclusion in \eqref{eq:subdiff} is nonempty and can be shown to be contained in the \subdiff{} of $g$ at $\bar{\Mvar}$, following the arguments in \citep[Example 3.1.1]{bertsekas2015convex}.
    To see why the inclusion in \eqref{eq:subdiff} holds, let $\bar{K}\in\KcalStar(\bar{\Mvar})$.
    Then, it follows that for any $\Mvar\in\Sset_+^\Nx\times\Sset_+^\Nx$,
    \begin{align*} 
        g(\Mvar) &= \min_{K\in\Lcal} \traceB{\big}{G(K)^\top \Mvar} 
        \le \traceB{\big}{G(\bar{K})^\top \Mvar}  
        = g(\bar{\Mvar}) + \traceB{\big}{G(\bar{K})^\top (\Mvar-\bar{\Mvar})},
    \end{align*}
    where the last equality follows from the definition of $g$ since $\bar{K}\in\KcalStar(\bar{\Mvar})$.
    Hence, $G(\bar{K})$ is a \subgrad{} of the function $g$ at $\bar{\Mvar}$, i.e., $G(\bar{K})\in\partial g(\bar{\Mvar})$.
   Since the \subdiff{} $\partial g(\bar{\Mvar})$ is a convex set in general (it is the intersection of closed halfspaces by definition), the convex hull of the \subgrads{} $G(\bar{K})$ for all $\bar{K}\in\KcalStar(\bar{\Mvar})$ must be contained in $\partial g(\bar{\Mvar})$, proving the inclusion in \eqref{eq:subdiff}.\\

    \noindent \emph{Proof of (iii):} Fix  $\bar{\Mvar} \in \Scal:=\{(\MvarA,\,\MvarB) \in\Sset_+^\Nx\times\Sset_+^\Nx \,\vert\, \MvarA+\MvarB\succ 0\}$ and let $$\phi(\Mvar,\,K):=\traceB{\big}{G(K)^\top \Mvar}  = \traceB{\big}{(\Mvar_1 + \Mvar_2)C(K)},$$
    where $\Mvar:=(\MvarA,\,\MvarB)$.\\[-.08in]

    \noindent \textbf{(Step 1) Uniqueness of the minimizer.} For any $\Mvar \in \Scal$, the function $K \mapsto \phi(\Mvar,\, K)$ is strongly convex on $\Rset^{m \times n}$, and therefore, its restriction to the linear subspace $\Lcal \subseteq \Rset^{m \times n}$ has a unique minimizer $\KcalStar(\Mvar) := \argmin_{K \in \Lcal} \phi(\Mvar, \, K)$.\\[-.08in]
    
    \noindent \textbf{(Step 2) Local compact reduction of the feasible set.}  Since $\Scal$ is an open set and eigenvalues are continuous functions, there is an open neighborhood $\Ncal$ of $\bar{\Mvar}$ and constants $\alpha, \beta > 0$ such that  for all $\Mvar\in \Ncal$, $$\alpha I \preceq \Mvar_1 + \Mvar_2 \preceq \beta I.$$
    For $\Mvar \in \Ncal$, we have
    \begin{align*}
        g(\Lambda) \leq \phi(\Mvar, \, 0) = \traceB{\big}{(\Mvar_1 + \Mvar_2)C(0)} \leq \beta \traceB{\big}{{\Kor}^\top D \Kor}.
    \end{align*}
    On the other hand, for any $\Mvar \in \Ncal$, we have 
    \begin{align*} 
        g(\Mvar) = \phi( \Mvar, \, \KcalStar(\Mvar)) \geq \alpha \lambda_{\rm min}(D) \| \KcalStar(\Mvar) - \Kor\|_2^2,
    \end{align*}
    where $\lambda_{\rm min}(D)>0$ denotes the smallest eigenvalue of the positive definite matrix $D$.
    Therefore, any minimizer $\KcalStar(\Mvar)$ (for $\Mvar \in \Ncal$) satisfies
    \begin{align*}
        \|\KcalStar(\Mvar) - \Kor\|_2 \leq R := \sqrt{\frac{\beta \traceB{\big}{K^{o \top} D \Kor}}{\alpha \lambda_{\rm min}(D)}} < \infty.
    \end{align*}
    Define the compact set $ \bar{\Lcal} :=  \{ K \in \Lcal \, | \,  \|K - \Kor\|_2 \leq R \}.$ Then, for all $\Mvar \in \Ncal$,
    \begin{align*}
        g(\Mvar) = \min_{K \in \bar{\Lcal}} \phi(\Mvar, \, K),
    \end{align*}
    since  every minimizer $\KcalStar(\Mvar)$ lies in the compact set $\bar{\Lcal}$, i.e., $ \{ \KcalStar(\Mvar) \, | \, \Mvar \in \Ncal\} \subseteq \bar{\Lcal}$.\\[-.08in]
    
    \noindent \textbf{(Step 3) Application of Danskin's theorem.}   The function $\phi(\Mvar,\,K)$ is jointly continuous in $(\Mvar,\,K)$ over the set $\Ncal \times \bar{\Lcal}$, and the map  $\Mvar \mapsto \phi(\Mvar,\, K)$ is linear (hence concave and differentiable) in $\Lambda$ for each $K \in \bar{\Lcal}$. Moreover, as shown in Step 1, the minimizer  $\KcalStar(\Mvar) := \argmin_{K \in \bar{\Lcal}} \phi(\Mvar, \, K)$ is unique for each $\Mvar \in \Ncal$. Since $\bar{\Lcal}$ is a compact set, Danskin's theorem \citep[Proposition B.25]{Bertsekas99} implies that 
    the function $g(\Mvar) = \min_{K \in \bar{\Lcal}} \phi(\Mvar, \, K)$ is differentiable on $\Ncal$, and for each $\Mvar \in \Ncal$,
    \begin{align*}
        \nabla g(\Mvar) = \nabla_{\Mvar} \phi(\Mvar, \, \KcalStar(\Mvar))=
        G(\KcalStar(\Mvar)),
    \end{align*}
    where $ \nabla_{\Mvar} \phi(\Mvar, \, K) $ denotes the gradient of $\phi$ with respect to its first argument $\Mvar$.
    In particular, $\nabla g(\bar{\Mvar}) = G(\KcalStar(\bar{\Mvar}))$. Since $\bar{\Mvar} \in \Scal$ was arbitrary, it follows that  $g$ is differentiable at any point in $\Scal$, concluding the proof.

\subsection{Proof of Theorem \ref{thm:pd_iterates}}

 We first prove the  inequality $\MvarB^\idx \succeq \SigmaNom$ by induction on the iteration  index $\idx$. The statement is true for the base case ($\idx = 0$), since $\MvarB^0 = \SigmaNom$ by definition. For the inductive step, suppose that $\MvarB^\idx \succeq \SigmaNom$ for some $\idx\ge 0$.  We need to show that $\MvarB^{\idx+1} \succeq \SigmaNom$. 
In Eq. \eqref{eq:proj_simp} in App. \ref{sec:proj}, it is shown that 
the projected \subgrad{} update rule for $\MvarB^{\idx+1}$ can be expressed as
\begin{align*}
    &\MvarB^{\idx+1} = \SigmaNom \, + \, \Pi_{\Scal_{r_2}^p}(\MvarB^\idx - \SigmaNom + \eta^\idx C(K^\idx) ),
\end{align*}
where $\Scal_{r_2}^p := \{X \in \Sset^{\Nx} \, | \, \| X \|_p \leq r_2\}$ is the Schatten $p$-norm ball of symmetric matrices centered at the origin. By the induction hypothesis, it holds that $\MvarB^\idx - \SigmaNom \succeq 0$, and since $\eta^\idx C(K^\idx) \succeq 0$, their sum remains positive semidefinite. 
Moreover, because the projection $\Pi_{\Scal_{r_2}^p}$ preserves positive semidefiniteness (see Eq. \eqref{eq:psd_proj} in App. \ref{sec:proj}), it follows that  $\MvarB^{\idx+1} \succeq \SigmaNom$, proving the desired claim.

Next, we prove that the function $g$ is differentiable on the sequence $\{\Mvar^\idx\}_{\idx=0}^\infty$ generated by the projected supergradient method \eqref{eq:subgrad_update_1}-\eqref{eq:subgrad_update_2}, when $\SigmaNom\succ 0$.
Recall from Theorem \ref{cor:max_sigma} that the dual function $g$ is differentiable at any point $(\MvarA,\,\MvarB)\in\Sset_+^{\Nx}\times \Sset_+^{\Nx}$ satisfying $\MvarA+\MvarB\succ 0$.
Using the previously established inequality $\MvarB^\idx \succeq \SigmaNom$ and the assumption that $\SigmaNom\succ 0$, it follows that $\MvarA^\idx+\MvarB^\idx \succeq \MvarA^\idx + \SigmaNom \succ 0$  for all $\idx\ge 0$, proving the desired result.

\section{Semidefinite Programming Reformulation} \label{sec:sdp}
Building on Theorem \ref{thm:main_result}, one can reformulate the \DRROC problem \eqref{eq:mroc} as a semidefinite program (SDP). 
This can be done by representing each term in the objective  in \eqref{eq:K_star} in epigraph form. The leading term can be expressed as 
\begin{align}
    \label{eq:expr1} &\traceB{\big}{\SigmaNom C(K)} = \min_{X\in\Sset^{\Nx}} \big\{ \trace{X} \, \big| \,  X  \succeq   \SigmaNom^\frac12 C(K)  \SigmaNom^\frac12 \big\},
\end{align}
where $\SigmaNom^\frac12 \in\Sset_+^\Nx$ denotes the (symmetric) positive semidefinite square root of the covariance matrix $\SigmaNom$. 
Using the Schur complement condition for positive semidefiniteness, the nonlinear matrix inequality  in \eqref{eq:expr1} can be reformulated as a linear matrix inequality (LMI):
\begin{align} \label{eq:SDP1}
\begin{bmatrix}
        X & \SigmaNom^\frac12 (K-\Kor)^\top\\
        (K-\Kor) \SigmaNom^\frac12 & D^{-1}
    \end{bmatrix}\succeq 0.
\end{align}
The regularization term involving the Schatten $q$-norm of the matrix $C(K)$ can also be represented in epigraph form as
\begin{align}
    \label{eq:expr2} &\|C(K)\|_q=\min_{\gamma\in\Rset, \, Y\in\Sset^\Nx}\big\{\gamma\,\big|\, \gamma \ge \|Y\|_q, \, Y \succeq C(K)\big\}
\end{align}
for $ q \in [1, \infty]$. For the spectral norm ($q = \infty$), the identity in \eqref{eq:expr2} simplifies to $\|C(K)\|_\infty = \min_{\gamma\in\Rset}\big\{\gamma\,\big|\, \gamma I \succeq C(K)\big\}$. Once again, we employ the Schur complement condition for positive semidefiniteness to reformulate the nonlinear matrix inequality  $ Y \succeq C(K)$ as the following LMI:
\begin{align}  \label{eq:SDP2}
\begin{bmatrix}
        Y & (K-\Kor)^\top\\
        K-\Kor & D^{-1}
    \end{bmatrix}\succeq 0.
\end{align}
The final constraint in the epigraph formulation \eqref{eq:expr2}, $ \gamma \ge \|Y\|_q$, is also  SDP representable. 
In particular, \citet[Prop. 4.2.1]{ben2001lectures} show that this constraint can be expressed as a set of LMIs whenever 
$q$ is rational. This result, together with the LMI reformulations in  \eqref{eq:SDP1} and \eqref{eq:SDP2}, yields an equivalent SDP reformulation of problem \eqref{eq:K_star}.

\section{Implementation Details of the Dual Projected \Subgrad{} Method} \label{sec:algo_details}

The dual projected \subgrad{} method presented in Sec. \ref{sec:opt_methods}, and summarized in Algorithm \ref{alg:pga}, involves two key steps at each iteration: (i) the evaluation of a dual \subgrad{} and (ii) the computation of a matrix projection onto the set $\Mball := \big\{ \Sigma\in\Sset_+^{\Nx} \,\big\vert\, \|\Sigma\|_1 \le r_1 \big\} \, \times \, \big\{ \Sigma\in\Sset_+^{\Nx} \,\big\vert\, \|\Sigma-\SigmaNom\|_p \le r_2 \big\}$.
We discuss how to efficiently carry out these computations in App. \ref{sec:subgrad_comp} and \ref{sec:proj}, respectively.
Additionally, in App. \ref{sec:stop}, we provide a stopping criterion for Algorithm \ref{alg:pga} using the relative duality gap of the primal-dual iterates generated by the proposed algorithm.

\begin{algorithm}
\caption{Dual Projected \Subgrad{} Algorithm}
\label{alg:pga}
\begin{algorithmic}[1]
\STATE{\textbf{Input:} Step size sequence $\{\eta^\idx\}_{\idx=0}^\infty$ and tolerance $\epsilon>0$}
\STATE{\textbf{Output:} $K^{\idx^\star}$ such that $f(K^{\idx^\star}) \le (1+\epsilon) f(\Kdrro)$}
\STATE Set $\Mvar^0 = (0,\,\SigmaNom)$\;
\STATE Compute $K^0 \in \argmin_{K\in\Lcal} \trace{G(K)^\top \Mvar^0}$\; \label{alg:K0}
\STATE Set $i=0$\;
\WHILE{$ \min_{0\le j,k \le \idx} \, \textrm{RelGap}(K^j,\Mvar^k)  > \epsilon$} \label{alg:stop}
    \STATE Compute $\Mvar^{\idx+1} = \Pi_{\Mball} (\Mvar^\idx + \eta^\idx G(K^\idx))$\; \label{alg:proj}
    \STATE Compute $K^{\idx+1} \in \argmin_{K\in\Lcal} \trace{G(K)^\top \Mvar^{\idx+1}}$\;  \label{alg:Kk}
    \STATE Update $\idx \gets \idx + 1$
\ENDWHILE
\STATE Set $\idx^\star \in \argmin_{j \le \idx} f(K^j)$\;
\STATE Return $K^{\idx^\star}$
\end{algorithmic}
\end{algorithm}

\subsection{Computation of Dual \Subgrad} \label{sec:subgrad_comp}
To evaluate a \subgrad{} $G(K^\idx) = (C(K^\idx), \, C(K^\idx)) $ at each iteration $i$, the equality-constrained convex quadratic program in Line \ref{alg:Kk} of Algorithm \ref{alg:pga} needs to be solved at every iteration. This can be accomplished by solving the linear system of equations derived from the Karush-Kuhn-Tucker (KKT) optimality conditions using direct factorization methods, or iterative methods such as the conjugate gradient method, among others \citep[Section 16]{nocedal1999numerical}.

\subsection{Computation of Matrix Projection onto Schatten \texorpdfstring{$p$-norm}{p-norm} Balls} \label{sec:proj}
We explain how to compute the matrix projection in Line \ref{alg:proj} of Algorithm \ref{alg:pga}. Because all Schatten norms are unitarily invariant, the orthogonal  projection of a matrix onto the  Schatten $p$-norm ball $\Scal_r^p := \{ X \in \Sset^{\Nx} \, | \, \| X \|_p \le r \}$  can be expressed in terms of the orthogonal projection of the corresponding vector of singular values onto the $\ell_p$-norm ball $\Bcal_r^p := \{ x \in \Rset^{\Nx} \,|\, \| x\|_p \le r\}.$ The following result, adapted from \citep[Theorem 7.18]{beck2017first}, makes this relationship precise.

\begin{thm} \rm \label{thm:proj}
    Let $X\in\Sset^{\Nx}$ and denote its eigendecomposition by $X=U \diag(\lambda) U^\top$, where $U$ is an orthonormal matrix. 
    Then 
    \begin{align} \label{eq:proj}
        \Pi_{\Scal_r^p} (X) = U \diag\big(\, \Pi_{\Bcal_r^p} (\lambda)\, \big) U^\top.
    \end{align}
\end{thm}
Here, $\diag(x) \in \Rset^{n\times n}$ denotes the diagonal matrix whose diagonal is the given vector $x\in\Rset^n$.
As an immediate consequence of this result, when the matrix being projected is positive semidefinite, the resulting projection remains positive semidefinite. This is due to the fact that the orthogonal projection of a nonnegative vector onto the $\ell_p$-norm ball $\Bcal_r^p$ is also nonnegative \citep[Lemma A.1]{bahmani2013unifying}.
It follows that
\begin{align}\label{eq:psd_proj}
X\succeq 0  \ \implies  \ \Pi_{\Scal_r^p \cap \Sset^\Nx_+} (X) = \Pi_{\Scal_r^p} (X) .
\end{align}
Using the above property,  together with the facts that $\eta^\idx C(K^\idx) \succeq 0$, $\MvarA^\idx \succeq 0$, and $\MvarB^\idx \succeq \SigmaNom$ (cf. Theorem~\ref{thm:pd_iterates}), the matrix projection in Line \ref{alg:proj} of Algorithm~\ref{alg:pga} can be simplified to
\begin{align} \label{eq:proj_simp}
    \Pi_{\Mball}(\Mvar^\idx + \eta^\idx G(K^\idx)) 
    &= \begin{bmatrix}
         \Pi_{\Scal_{r_1}^1}(\MvarA^\idx + \eta^\idx C(K^\idx))\\
         \SigmaNom+\Pi_{\Scal_{r_2}^p}(\MvarB^\idx + \eta^\idx C(K^\idx)-\SigmaNom)
    \end{bmatrix}.
\end{align}
To compute the matrix projections in the right-hand side of Eq. \eqref{eq:proj_simp}, one first needs to compute eigendecompositions of the matrices appearing inside  each of the projection operators. Theorem~\ref{thm:proj} can then be applied to express each matrix projection in terms of a vector projection of the resulting eigenvalues onto the appropriate $\ell_p$-norm ball. Efficient algorithms  for computing these vector projections are detailed in \citep{won2023unified}. Note that the overall computational complexity of the projection step is dominated by the eigendecompositions, which scales as $\mathcal{O}(\Nx^3)$ for dense matrices.

\subsection{Stopping Criterion} \label{sec:stop}
Here, we introduce a stopping criterion for the sequence of iterates generated by Algorithm \ref{alg:pga}.
Given a feasible dual point $\Mvar \in \Mball$, 
one can bound the suboptimality of any feasible control policy $K \in \Lcal$ using the \emph{relative duality gap}, defined as
\begin{align} \label{eq:relgap}
    \textrm{RelGap}(K,\Mvar) := (f(K)-g(\Mvar))/g(\Mvar).
\end{align}
For the relative duality gap \eqref{eq:relgap} to be well defined, we require that $g(\Mvar) > 0$. It can be shown that this condition holds whenever $\Kor\not\in\Lcal$ and $\MvarA+\MvarB\succ 0$.  As noted in Sec. \ref{sec:pga}, if Algorithm \ref{alg:pga} is initialized at $\Mvar^0 = (0,\,\SigmaNom)$ with  $\SigmaNom\succ 0$, then all subsequent iterates $\{\Mvar^\idx=(\MvarA^\idx,\,\MvarB^\idx)\}_{i=1}^\infty$ satisfy $\MvarA^\idx + \MvarB^\idx \succ 0$.
By weak duality, \eqref{eq:relgap} provides an upper bound on the \emph{relative optimality gap} of the control policy $K$, i.e., 
\begin{align} \label{eq:upper_bd}
    \textrm{RelGap}(K,\Mvar) \ge \frac{f(K) - \min_{K'\in\Lcal}f(K')}{\min_{K'\in\Lcal} f(K')}.
\end{align}
Therefore, one can stop the algorithm when this upper bound falls below a desired tolerance. 
Since the proposed algorithm is not guaranteed to produce monotonic sequences of primal and dual objective values, we compute  the relative duality gap at each iteration $i$  (Line \ref{alg:stop} of Algorithm \ref{alg:pga}) using the best primal-dual iterates computed up to that point: 
\begin{align} \label{eq:bestGap}
    \min_{0\le j,k \le \idx} \textrm{RelGap}(K^j,\Mvar^k).
\end{align}

\section{Additional Numerical Experiments} \label{sec:additional_exp}

In App. \ref{sec:conv_scalability}, we examine the convergence and scalability of the proposed algorithm, and compare its computation time against that of a state-of-the-art interior-point solver applied to the SDP reformulation given in App. \ref{sec:sdp}. In App. \ref{sec:add_exp_with_mean}, we provide  additional numerical experiments for a variant of the system considered in Sec. \ref{sec:experiments} in which the disturbance mean is modeled as nonzero and uncertain.

\subsection{Convergence and Scalability Experiments}  \label{sec:conv_scalability}

\begin{figure*}
    \centering
    \subfigure[Convergence of Algorithm \ref{alg:pga} \label{fig:relGap}]{\includegraphics[width=0.48\columnwidth]{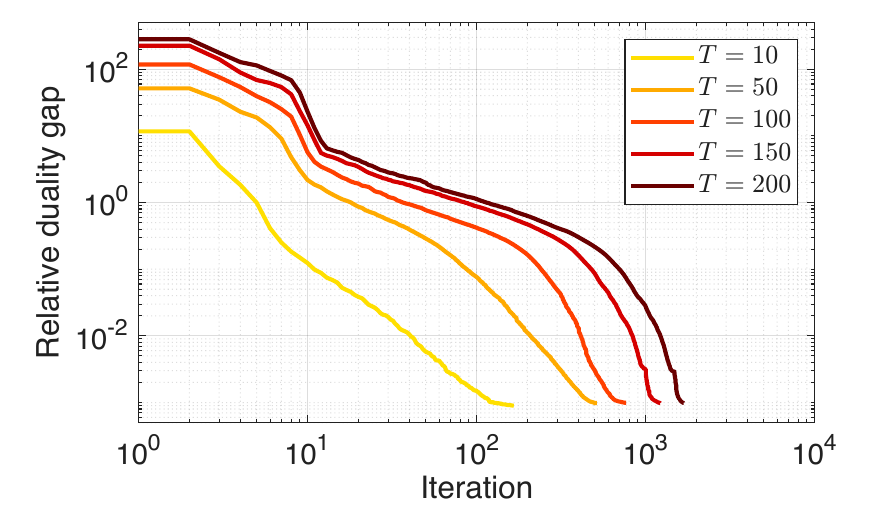}}
    \subfigure[Computation time comparison \label{fig:compTimes}]{\includegraphics[width=0.51\columnwidth]{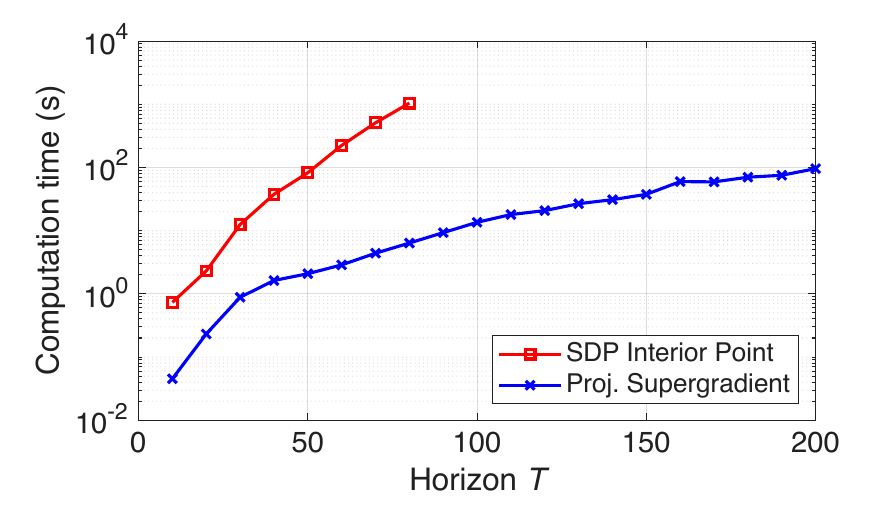}}
    \caption{(a) Relative duality gap \eqref{eq:bestGap} (averaged over ten trials) versus the iteration count for the dual projected \subgrad{} method (Algorithm \ref{alg:pga}) for different control horizons $T$. (b) Total execution time (averaged over ten trials) as a function of the control horizon $T$ for the SDP interior point method (red line, square markers) and the dual projected \subgrad{} method (blue line, cross markers).} 
\end{figure*}

We examine the convergence and scalability of the proposed dual projected \subgrad{} method (Algorithm \ref{alg:pga}) when the covariance uncertainty set is defined in terms of the nuclear norm ($p=1$).
The numerical experiments are carried out using the system described in Sec. \ref{sec:experiments} for a correlation parameter value of $\rho=0$.
Algorithm \ref{alg:pga} is implemented using the adaptive step size rule:
$$\eta^\idx = \min \left(1,\,\frac{\sqrt{2} \, \|\Mvar^\idx-\Mvar^{\idx-1}\|_2}{ \|\nabla g(\Mvar^\idx) - \nabla g(\Mvar^{\idx-1})\|_2} \right).$$
The inverse step size $1/\eta^\idx$ can be interpreted as a local estimate of the  Lipschitz constant of the gradient $\nabla g$ on $\Sset_+^\Nx \times \Sset_+^\Nx$, truncated from above at one.\footnote{
This adaptive step size rule is closely related to the Barzilai-Borwein step sizes used in spectral projected gradient methods \citep{birgin2000nonmonotone}.
Gradient methods utilizing such adaptive step sizes enjoy  convergence guarantees when combined with saturation limits and nonmonotone line search schemes \citep{wang2005convergence}.}
We demonstrate empirically that the proposed algorithm converges when using this step size sequence for different control horizons $T$. We also compare the computation time of solving problem \eqref{eq:K_star} using this method with that of an interior point optimizer applied to the SDP reformulation described in App. \ref{sec:sdp}. 

In our experiments, the ambiguity set radius is scaled with the control horizon as $r=T$, and we vary the control horizon from $T = 10$ to $T = 200$ by increments of $10$. 
For each horizon~$T$, we conduct ten independent trials, where we draw a new training sample every trial, and solve problem \eqref{eq:K_star} using the two optimization methods under comparison.
We record the time required to solve the SDP reformulation of problem \eqref{eq:K_star}, which we solve in MATLAB using CVX (version 2.2) and the MOSEK interior point solver (version 9.1.9) with a relative gap tolerance of $10^{-3}$ and infeasibility tolerance of $10^{-8}$.
For the proposed \subgrad{} method, we record the time required to satisfy a relative gap tolerance of $\epsilon=10^{-3}$. 
The orthogonal projection onto the $\ell_1$-ball (which is needed to compute the matrix projection in Line \ref{alg:proj} of Algorithm \ref{alg:pga} via the procedure described in App. \ref{sec:proj}) is implemented using the algorithm of \citet{duchi2008efficient}.
A laptop with an Intel Core i7-8665U processor and 16 GB of RAM was used for all computations.

We first examine the convergence of the projected \subgrad{} method.
Figure \ref{fig:relGap} depicts the relative duality gap according to the bound in \eqref{eq:bestGap} (averaged over  ten trials) versus the iteration count for control horizons $T=10,\,50,\,100,\,150,\,\text{and }200$.
While the number of iterations needed for the relative duality gap to fall below $\epsilon=10^{-3}$ increases with the control horizon $T$, the algorithm consistently converges within a moderate number of iterations  for all  values of $T$ considered.

Next, we evaluate the scalability of the projected \subgrad{} method relative to the SDP interior-point solver.
Figure \ref{fig:compTimes} reports the average computation times (over ten trials) for solving problem \eqref{eq:K_star} as a function of the control horizon $T$.
We exclude the SDP interior-point solver times for $T>80$, because the SDP  becomes too large for the solver to handle within the available memory.
Across all horizons where the comparison is feasible, the projected \subgrad{} method attains solutions within the prescribed tolerance substantially faster than the interior-point method.
For example, at $T=80$, the interior-point method takes over 16 minutes on average, whereas the projected \subgrad{} method converges in under eight seconds on average.

\subsection{Additional Experiments with Disturbance Mean Uncertainty} \label{sec:add_exp_with_mean}

We report additional numerical results for a modified version of the system in Sec. \ref{sec:experiments}, where the ground-truth disturbance mean is nonzero and uncertain. Specifically, we consider a  disturbance process of the form 
$$\widetilde{w} =  w  \, + \,  \Big(\frac{1}{10} \Big) \mathds{1},$$
where  $w$ is the same AR(1) disturbance process specified in Sec. \ref{sec:experiments} and $\mathds{1} \in \Rset^\Nx$ is the vector of all ones.
Given a finite training sample of $N=\Nx+1$ i.i.d. disturbance trajectories $\widetilde{w}^{(1)},\dots,\widetilde{w}^{(N)}\in\Rset^\Nx$,  we compute the nominal mean and covariance matrix via the  unbiased empirical estimators:
\begin{align*}
    \muNom = \frac1N \sum_{i=1}^N \widetilde{w}^{(i)}, \qquad \SigmaNom = \frac{1}{N-1}\sum_{i=1}^N (\widetilde{w}^{(i)}-\muNom)(\widetilde{w}^{(i)}-\muNom)^\top.
\end{align*}
Using this setup, we repeat the numerical experiments of Sec. \ref{sec:experiments} to evaluate how the out-of-sample performance varies with the mean and covariance ambiguity set radii $(r_1, \, r_2)$.
For the controllers proposed in this paper, we vary the ambiguity set radii $(r_1,\,r_2)$ over a finite grid on $ [10^{-4},\,10^4]\times [10^{-4},\,10^4]$, and recompute the controllers for every grid point. 

\begin{figure}
    \centering
    \includegraphics[width=\linewidth,trim=0cm 40.15cm 0cm 0cm,clip]{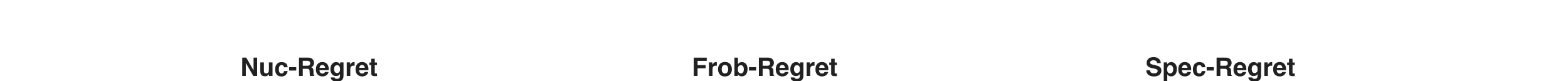}
    \includegraphics[width=\linewidth,trim=0cm 1.25cm 0cm 1.5cm,clip]{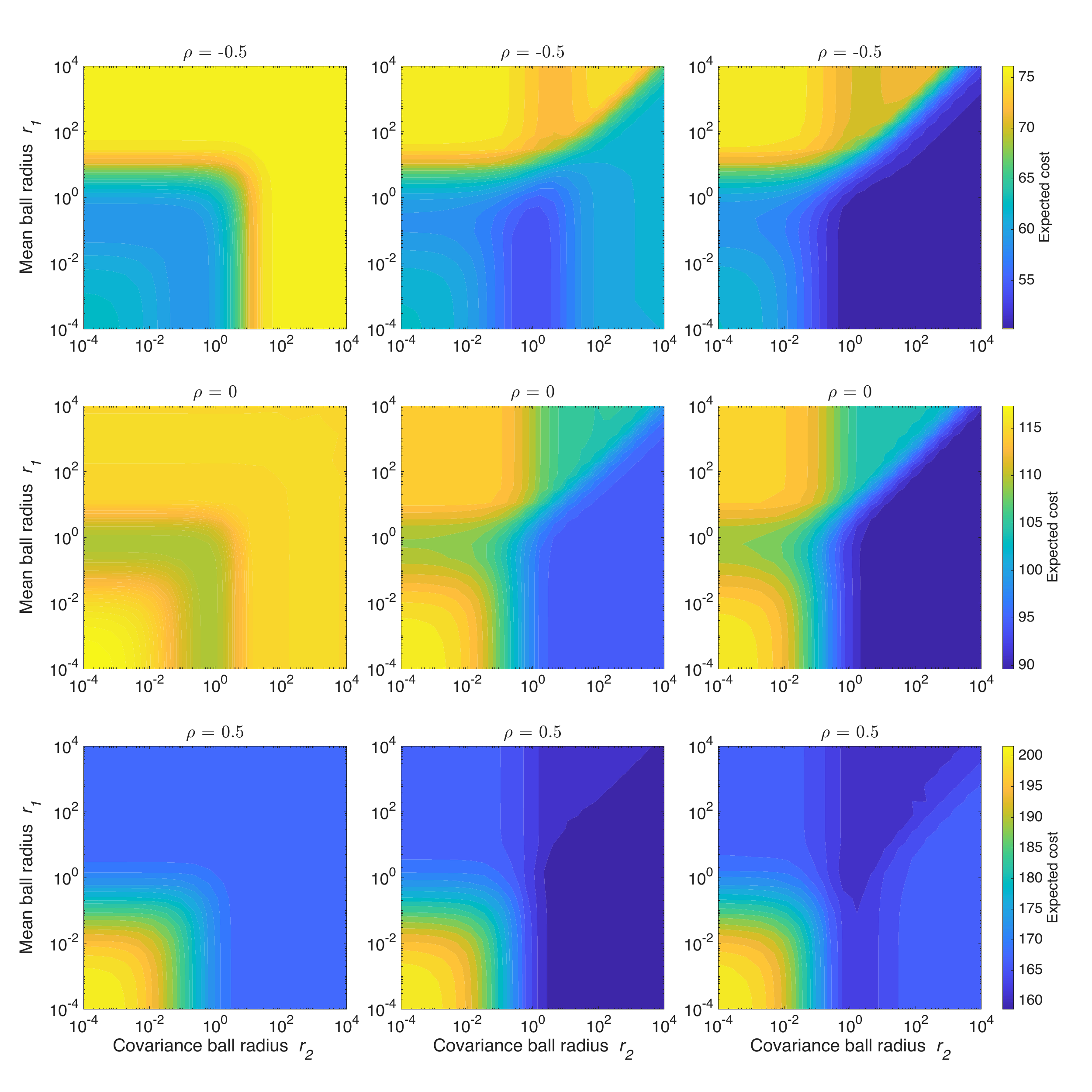}
    \caption{Expected cost under the ground-truth distribution  as a function of the ambiguity set radii $(r_1,\,r_2)$ for the \textsc{Nuc-Regret} (left column), \textsc{Frob-Regret} (middle column), and \textsc{Spec-Regret} (right column) controllers. Across rows, $\rho$ varies over $\{-0.5,\, 0, \, 0.5\}$. The expected costs are averaged over 50 independent trials, where, in each trial, the nominal mean and covariance are estimated from a new training sample drawn from the ground-truth distribution.}
    \label{fig:heatmap}
\end{figure}

Fig. \ref{fig:heatmap} shows the expected cost of the resulting controllers under the ground-truth distribution (averaged over 50 independent trials) as a function of the ambiguity set radii $(r_1,\,r_2)$ for three different correlation coefficients $\rho \in \{-0.5,\,0,\,0.5\}$. As expected, the out-of-sample performance is sensitive to both radii, highlighting the need to  jointly calibrate $(r_1, \, r_2)$ to achieve the best out-of-sample performance.

\begin{figure}
    \centering
    \includegraphics[width=1.1\linewidth,trim=1.75cm 0.5cm 0cm 0.45cm,clip]{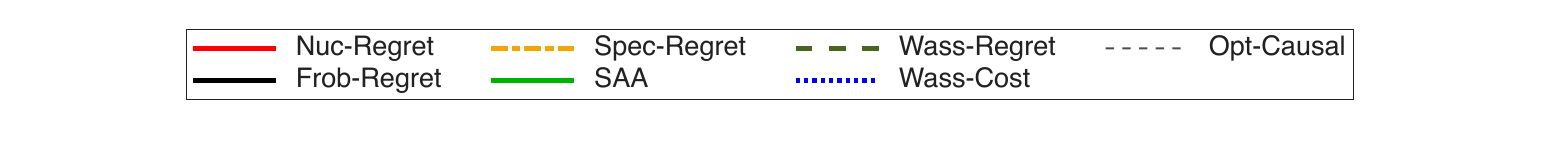}
    \includegraphics[width=\linewidth,trim=3.9cm 0cm 4cm 0cm,clip]{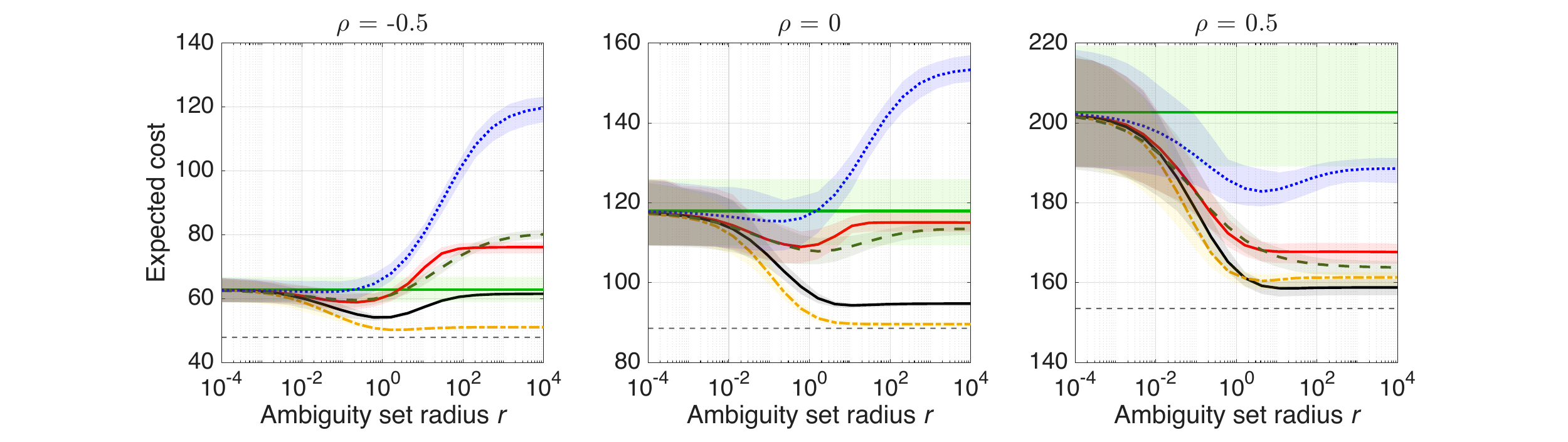}
    \caption{Expected controller cost under the ground-truth distribution versus the ambiguity radius $r$, for  $\rho \in \{-0.5,\, 0, \, 0.5\}$. 
    The Wasserstein-based controllers are computed with radius $\sqrt{r}$, whereas the remaining distributionally robust controllers are computed for each $r$ at the  pair $(r_1,r_2)$ satisfying $r_1+r_2=r$ that minimizes the ground-truth  expected cost.
     Expected costs are averaged over 50 independent trials, where, in each trial, the nominal mean and covariance are estimated from a new training sample drawn from the ground-truth distribution.
    The shaded regions depict the empirical 20\textsuperscript{th}--80\textsuperscript{th} percentile range.
    }
    \label{fig:ecost_vs_r_nzmean}
\end{figure}

Fig.  \ref{fig:ecost_vs_r_nzmean} compares the expected cost under the ground-truth distribution of the proposed controllers with that of the Wasserstein-based controllers that are parameterized by a single radius $\sqrt{r}$. To evaluate the proposed controllers on the same one-dimensional scale, we sweep over $r\in[10^{-4},\,10^4]$ and, for each value of $r$, we select the pair $(r_1, \, r_2)$  satisfying $r_1+r_2 = r$ that minimizes the  expected cost under the ground truth distribution. For the three different correlation coefficients considered, the \textsc{Spec-Regret} and \textsc{Frob-Regret} controllers proposed in this paper achieve the best out-of-sample performance.

